\documentclass[12pt,a4paper]{article}
\usepackage{graphicx} 
\usepackage{amsmath,amssymb,mathtools,geometry}
\usepackage{amsthm}

\usepackage{enumitem}
\usepackage{booktabs}
\usepackage[colorlinks=true, citecolor={blue}, linkcolor={blue}, urlcolor={blue}]{hyperref}
\usepackage{xcolor}
\usepackage[round,authoryear]{natbib}
\usepackage{float}
\usepackage{amsmath,amssymb}
\usepackage{algorithm}
\usepackage{algpseudocode}
\algrenewcommand{\algorithmiccomment}[1]{\hskip1em$\triangleright$ #1}
\usepackage{longtable,array}
\newcolumntype{R}[1]{>{\raggedleft\arraybackslash}p{#1}}
\usepackage{tabularx}
\usepackage{array}
\newcolumntype{R}[1]{>{\raggedleft\arraybackslash}p{#1}}
\newcolumntype{Y}{>{\raggedleft\arraybackslash}X}
\usepackage{tikz}
\usetikzlibrary{arrows.meta, positioning, fit, backgrounds, calc}

\usepackage{authblk}

\title{Designing an Adaptive Robust and Evolutionary-Game Framework for Minimum Shipment Quotas in Capacity-Constrained Freight Markets}

\author[a]{Hossein Mirzaee\thanks{Corresponding Author}}
\author[b]{Ali Keyvandarian}
\author[c]{Hamed Samarghandi}
\affil[a]{\small{University of British Columbia, Vancouver, British Columbia, Canada, V6T 1Z4}}
\affil[b]{\small{Department of Industrial Engineering, Dalhousie University, Halifax, NS, Canada, B3H 4R2}}
\affil[c]{\small{Edwards School of Business, University of Saskatchewan, Saskatoon, SK, Canada, S7N 5A7}}
\date{}

\usepackage{setspace}
\begin{document}

\maketitle

\begin{abstract}
Minimum shipment quota regulations ensure that socially critical commodities are not excluded from capacity-constrained freight markets by profit-driven allocation. Existing models of such regulation are deterministic and static, assuming perfect market knowledge and instantaneous convergence to a rational equilibrium. This paper relaxes both assumptions through a framework combining adaptive robust optimization, multi-objective optimization, and evolutionary game theory. A robust optimization model first characterizes the regulator's quota-setting problem under uncertain commodity availability, trading social welfare against enforcement cost to generate a Pareto frontier of quota policies. This frontier then defines the regulator's strategy space in a three-population evolutionary game in which the regulator, suppliers, and carrier co-evolve under replicator dynamics, so that the prevailing policy emerges from the interaction of the three players. The framework is calibrated to a Canadian rail freight case study that covers crude oil, corn, barley, and oat. Quota regulation raises social welfare from 65.0 to 288.2 million CAD, at an enforcement cost of about 2\% of the gain. Two findings carry particular policy weight: capacity shortfalls fall entirely on unprotected commodities the carrier ships voluntarily, sparing both the highest-margin and the protected products; and a fall in transport cost can end the regulation while the exclusion it was designed to correct persists, because the margin observable before bidding diverges from the margin prevailing once suppliers adapt.
\end{abstract}


\section{Introduction}\label{sec:introduction}

Railway freight transportation is considered the main backbone of exporting countries, as they are responsible for the transportation of bulk goods to ports of consumption, being one of the main types of trade for countries rich in grain, mineral and energy resources \citep{aar2025}. In Canada, where there is a large distance between the producing regions and  demand points, railway carriers transport more than 325 million tons of cargo per year, being the major type of transport for this country \citep{statcan2025}. However, the transportation capacity is determined by the available number of locomotives, cars, tracks, and terminals for loading, which acts as a scarce and shared resource that is competitively assigned to different commodity flows. If the capacity is less than the demand, then the market mechanism of allocation will provide priority to the commodities with the highest transport margins, leaving important products without any opportunity to reach markets \citep{taherkhani2022}. This situation creates the conflict between profit maximization and social welfare, leading to regulatory interventions through tools such as minimum shipment quotas \citep{pazoki2024}.

This conflict was dramatically illustrated in Canada during the 2013-2014 crop year, when a record grain harvest coincided with a rapid increase in crude oil shipments by rail due to limited pipeline capacity and high oil prices, while unusual winter caused rail shipments to decrease even further \citep{cta2014}. Since oil producers were willing to pay much higher prices for rail services, rail companies preferred to use the remaining capacity for oil transportation, leaving an estimated 18 million tons of grains in elevators and on farms, thus causing losses for Canadian farmers between two and four billion Canadian dollars (CAD) \citep{nytimes2014}. After this event, the Federal government passed the Fair Rail for Grain Farmers Act in May 2014, requiring the two largest rail companies in Canada to deliver at least 500,000 tons of grain per week and introducing financial penalties for non-compliance \citep{canada2014}. However, stakeholders deemed the regulation unacceptable, as grain producers believed that the quotas were too low, and oil producers and carriers did not want government interference in market allocation \citep{pazoki2024}. This paper formulates the problem of designing a minimum quota regulation in terms of economic efficiency and political viability in a competitive environment.

In the rail freight regulation problem, there are three strategic players with completely opposite objectives. The regulator tries to maximize social welfare by ensuring the delivery of certain goods to markets despite their exclusion by market dynamics due to low margins as these goods are critical for the stability of the food, energy or rural economy \citep{demirci2017}. Suppliers, producers of grains, oil and other commodities compete with each other to offer maximum prices for transportation to secure enough capacity for shipment, acting as Bertrand competitors with capacity constraints \citep{pazoki2024}. The rail carrier works as a profit maximizing follower allocating capacity to the most profitable flows after setting quotas. The problem cannot be formulated as a one-optimizer programming model since the players do not control the outcome of the problem; therefore, this problem should be studied as a strategic interaction game \citep{adler2021}.

The real world rail freight market is characterized by great uncertainty, elevating the degree to which classical optimization approaches can make serious mistakes. Commodity demand is affected by harvest yield, global shocks, and export conditions \citep{roberts2013}. The cost of production fluctuates depending on the price of fuel, labor rates, and input costs \citep{hamilton2009}. Rail capacity fluctuates seasonally due to weather, equipment availability, and infrastructure maintenance \citep{taherkhani2022}. Large scale disruptions, from pandemics to weather, lead to additional uncertainty and motivate research devoted to strategies of supply chain resilience under stochastic conditions \citep{mirzaee2024comparing}. Even the enforcement of the policy itself is associated with uncertainties, caused by monitoring efforts, political pressure, and costs of compliance, which evolve over time. Neglecting this fact leads to the fragility of policies which work fine in nominal conditions but badly in adverse situations, where intervention is called for. The robust design of regulation requires approaches that take into account uncertainty in the parameters \citep{bertsimas2004}.

Various approaches to decision making under uncertainty have been developed; for instance, stochastic programming, chance-constrained optimization, and robust optimization \citep{gabrel2014}. Stochastic programming requires a probability distribution for uncertain parameters that is usually unknown for policy design problems when there is little historical data or when the system undergoes structural changes \citep{shapiro2014}. Chance-constrained approaches allow for the violation of constraints with some probability, but require distributional assumptions and are sometimes computationally hard for large-scale problems \citep{charnes1959}. On the other hand, the robust approach only requires specification of the uncertainty set in which the parameters may fluctuate and seeks to make effective decisions under the worst realization in this set \citep{soyster1973, mirzaee2023}. The worst-case approach suits regulators, who must consider outcomes under adverse conditions and protect society from bad scenarios \citep{hansen2008}. In this work, we adopt the adaptive variant of robust optimization (ARO), proposed by \citet{bental2004}, which facilitates the two stage decision making process of the study, where the regulator quotas are set in advance, while the carrier capacity allocation is a recourse decision made after the resolution of uncertainty. This two-stage structure is perfectly captured by ARO.

Game-theoretic models of rail freight regulation, especially the recent Stackelberg model presented by \cite{pazoki2024} facilitate understanding the equilibrium behavior of suppliers and the resulting socially optimal quota policy. However, such models have one common shortcoming: they assume that all players are fully rational, have complete information about payoffs, and can converge to the Nash equilibrium instantaneously. Whereas in the real world freight market, suppliers do not have complete information about costs, valuations, and competitor strategies. Instead, suppliers observe the result in the market over some periods of time and slowly learn how to behave properly \citep{fudenberg1998}. The carrier does the same with respect to the allocation rules. Consequently, the regulator adjusts their position in response to market behavior. The One-step static equilibrium approach cannot model this learning process or analyze the stability of candidate equilibria, which endure in such a learning process.

Evolutionary game theory (EGT) is a method originally developed in biology and is now widely used in economics and operations research \citep{weibull1995}. In this framework, strategies do not represent the choice of rational individuals, but the prevailing behavior in the population. In EGT, the share of each strategy changes in time according to replicator dynamics, in which the share of strategies leading to above-average payoffs spread in population, while the others shrink. EGT allows the modeler to incorporate bounded rationality and the stability of equilibria, and accepts mixed equilibria and limit cycles as natural possible outcomes \citep{sandholm2010}. Compared to pure static equilibrium approaches, embedding the rail freight regulation problem with three players in the EGT framework with co-evolution of all players' populations leads to the development of a model that represents policy design in the real world more accurately.

In light of the issues discussed above, this paper examines the following three inter-related research questions:
\begin{enumerate}
    \item How can a regulator determine minimum shipment quotas to maximize social welfare in the presence of uncertainties in demand for commodities, transportation costs, and other variables?
    \item Within the group of Pareto-efficient robust quota rules, what is the rule that is expected to sustain evolutionary competition among suppliers and carriers while they adapt their behaviors?
    \item What is the sensitivity of the equilibrium rule to changes in economic parameters such as the relative social value of commodities, enforcement cost, uncertainty in demand, rail capacity, and market prices of commodities?
\end{enumerate}

The contributions of this article are as follows.
\begin{enumerate}
    \item Two-stage robust optimization, multi-objective optimization, and evolutionary game theory are integrated into a single framework for designing quota regulation in rail freight. The robust stage generates Pareto-efficient quota policies that trade social welfare against enforcement cost. The evolutionary stage then selects the policy that persists as suppliers and the carrier adapt their strategies.

    \item In existing evolutionary models of regulation, the regulator chooses from a fixed set of policy levers. Here, the regulator's strategies are the policies on the robust Pareto frontier, so each policy in the evolutionary game is robust to uncertainty and not dominated by any other.

        \item The model exposes a form of regulatory failure. If protection is granted on the basis of margins observed before bidding, a commodity can remain excluded at the bids suppliers actually post. In the case study this happens twice: a fall in transport cost ends the regulation of oat, and a higher social weight on barley never leads to its protection.
\end{enumerate}

When applied to the Canadian case, the framework generates three policy implications, as follows.
\begin{enumerate}
    \item Quota regulation raises social welfare from $65.0$ million CAD to $288.2$ million CAD at an enforcement cost of $4.8$ million CAD, about 2\% of the gain. The carrier bears most of the private cost, and the oat supplier, excluded without regulation, gains the most.

    \item Capacity is the parameter with the largest effect. When capacity contracts, the shortfall falls on the commodity that is neither protected nor the most profitable, which in this case is corn. 
    
    \item A fall in the carrier's transport cost can end the regulation while the exclusion it corrects persists. Lower costs make the maximum-bid margin of oat positive before bidding, so the robust stage sees no need for a quota, yet at equilibrium oat still bids zero and is not shipped.
\end{enumerate}

The remainder of this paper is organized as follows. Section~\ref{sec:literature} reviews the related literature. Section~\ref{sec:problem} describes the problem and states the modeling assumptions. Section~\ref{sec:solution-framework} presents the solution framework and the proposed algorithm. Sections~\ref{sec:aro} and~\ref{sec:egt} develop the adaptive robust optimization and evolutionary game-theoretic models, respectively. Section~\ref{sec:case} describes the case study, and Section~\ref{sec:results} reports the numerical results, as well as sensitivity analyses and managerial insights. Section~\ref{sec:conclusion} concludes the paper.

\section{Literature Review}\label{sec:literature}

This paper draws on five streams of literature: (i) transportation pricing and capacity allocation, (ii) game-theoretic models of freight markets, (iii) transportation regulation, (iv) robust optimization under uncertainty, and (v) evolutionary game theory in operations research.

\subsection{Transportation Pricing and Capacity Allocation}

The problems of pricing transportation services and allocating scarce transport capacity have been approached with a wide range of methodologies. Early and ongoing work is based on deterministic mathematical programming. \citet{azadian2018} developed a mixed-integer non-linear program to determine transportation prices based on the geographical grouping of service locations, while \citet{saif2021new} formulated a mixed-integer nonlinear supply chain network design model that captures transportation and capacity economies of scale. \citet{chang2018} and \citet{gu2018} formulated road pricing models that internalize congestion and environmental externalities through distance- and time-based tolls. 

Dynamic programming and revenue management form a second methodological family, particularly for carriers that must ration capacity across heterogeneous customer classes. \citet{amaruchkul2011} designed a dynamic programming-based revenue management scheme for air-cargo capacity allocation, \citet{xu2022} proposed revenue management schemes for high-speed rail systems shared between passengers and freight under deterministic and stochastic demand, and \citet{taherkhani2022} developed a tactical capacity planning model for a multi-stakeholder freight system with revenue management, analyzing how parameter changes propagate through the transportation plan. More recently, \citet{bao2026slot} formulated slot allocation across a regional multi-airport system as a tri-objective problem trading total schedule displacement against fairness to airlines and to airports, and generated the Pareto frontier by an $\varepsilon$-constraint procedure. 

Simulation-based approaches have also been employed where analytical tractability is limited. \citet{abril2008assessment} reviewed the main methods for railway capacity analysis, simulation among them, and presented an automated tool that runs several of these analyses. \citet{dingler2009impact} used train dispatching simulation on a single-track North American line and found that mixing intermodal, unit, manifest, and passenger trains causes more delay than homogeneous traffic, which reduces the capacity the line can offer. \citet{bablinski2016} built a simulation game based on real data from the Brighton Main Line to study capacity allocation between competing freight and passenger operators.

Although these methods yield valuable prescriptions for a single decision-maker, they share a structural limitation. They treat the focal decision-maker environment as exogenous. When several self-interested actors, such as producers, carriers, and regulators, make interdependent decisions, the outcome of each player's optimization depends on the anticipated behavior of the others. Game theory addresses precisely this interdependence, which explains its growing prominence as one of the most suitable methodologies for transportation pricing and capacity allocation problems \citep{adler2021}.

\subsection{Game-Theoretic Models of Freight Markets}

Game-theoretic treatments of freight markets initially focused on two-player interactions. \citet{zhou2009} studied pricing and competition between transportation service providers with empty equipment repositioning. \citet{toptal2011} examined pricing in a market where a truckload and a less-than-truckload carrier compete and demonstrated gains from cooperation. \citet{mozafari2015} modeled dynamic pricing among freight carriers competing in an oligopolistic network as a discrete-time dynamic Nash equilibrium. \citet{peng2016} applied a Gale--Shapley-based price game mechanism to a stable vessel-cargo matching in dry bulk shipping. 

Stackelberg leader--follower structures have proven especially natural for freight settings in which one party moves first. \citet{tamannaei2021} formulated a Stackelberg competition between a direct transportation system and intermodal service providers to find equilibrium prices in a multi-level freight market. However, many real transportation problems involve \emph{three} distinct classes of players whose objectives conflict, typically a public authority, one or more producers or shippers, and a transport service provider. A growing stream of research addresses such tri-level interactions. \citet{mirzaee2022} developed a three-player stochastic game between a manufacturer, a third-party emission verifier, showing how the government can set penalties to maximize social welfare while deterring collusion. Most directly related to this article, \citet{pazoki2024} modeled the Canadian rail freight market as a Stackelberg game in which suppliers compete \`a la Bertrand for the capacity of a profit-maximizing carrier, and characterized the regulator's optimal minimum-quota policy through a social utility maximization model. Our paper builds on this tradition of three-players but extends it in three directions that the static, deterministic literature does not address: 1)~parameter uncertainty, 2)~multiple regulatory objectives, and 3)~adaptive (rather than perfectly rational) behavior.

\subsection{Transportation Regulation}

The regulation of freight transportation has been studied through several instruments other than shipment quotas. The Hours-of-service rules were analyzed by \citet{goel2014}, who compared the safety and cost implications of the regulations in the European Union, the United States, and Canada. The Sustainability-oriented regulations on urban freight were reviewed by \citet{quak2008}. \citet{wei2020} examined route optimization under the low-carbon emission regulation. \citet{tiwari2021} studied freight consolidation and containerization under the carbon tax regulation. In a congested crowd sourced delivery market, \citet{li2026regulating} showed that no single regulatory instrument restores first-best efficiency, and that a commission cap must be paired with a bundling constraint to align the platform's decisions with the social optimum. \citet{demirci2017} investigated interventions designed to steer the demand towards socially preferred goods in the context of the California electric vehicle market. 

Closer to our setting, price-cap regulation, for which Canada's Maximum Revenue Entitlement is an example, has a long tradition in the economics of network industries \citep{beesley1989}. \citet{kang2000} showed that whether consumers benefit from tightening a price cap depends on demand interdependence. In contrast to these instruments, \emph{the minimum quota} regulation of transportation capacity has received little attention. \citet{pazoki2024} provided the first game-theoretic analysis of the minimum quota regulation in rail freight, motivated by Canada's \emph{Fair Rail for Grain Farmers Act}, and derived conditions under which a product deserves a quota. However, their analysis was deterministic and demand, prices, and costs are assumed to be static and known in the sense that all players reach equilibrium instantaneously and enforcement is costless. Both simplifications matter in practice, since quota policies must be set before harvest volumes, prices, and capacities are realized, and enforcing quotas consumes real administrative and political resources.

\subsection{Robust optimization in Supply Chains and Transportation}

In scenarios where the distributional information of the uncertainties is unknown or inaccurate, Robust optimization (RO) has been identified as the prevailing approach to decision-making. Following the foundational contributions of \citet{soyster1973}, \citet{bental2004}, and \citet{bertsimas2004}, RO seeks decisions that remain feasible and perform well for all realizations within a prescribed uncertainty set, and its adjustable (two-stage) variant distinguishes here-and-now decisions from recourse decisions made after uncertainty resolves \citep{bental2004}. RO applications in logistics and transportation have expanded rapidly. In freight transportation, \citet{zhang2026sharing} formulated a multi-carrier robust service network design problem in which carriers share orders or vehicles under demand uncertainty. \citet{guo2026integrated} developed a robust location-routing framework for rural passenger-freight transport under leasing-cost uncertainty. In the resilience literature, \citet{liu2024} designed robust actions to improve supply chain viability under disruption. \citet{zhou2024} analyzed the resilience of the integrated China--Europe freight transport network under heterogeneous demands. Recent transportation studies have also applied RO to intermodal routing under uncertain transit times and carbon constraints \citep{xu2025eurasian}. In the context of aviation decarbonization, \citet{keyvandarian2026robust} developed a robust adaptive optimization framework for the distribution of sustainable aviation fuel (SAF) in Canada under uncertain production scenarios. 

Collectively, the literature demonstrates that robust methods produce policies that degrade gracefully under adverse realizations, i.e., the property a politically accountable regulator demands. However, to the best of our knowledge, RO has not been applied to the design of transportation quota regulation, nor has any robust framework been coupled with a game-theoretic model of the regulated market's response.

\subsection{Evolutionary Game Theory in Operations and Transportation}

Evolutionary game theory (EGT), originating with \citet{smith1973} and systematized by \citet{weibull1995} and \citet{sandholm2010}, replaces the classical assumption of perfectly rational players with populations of boundedly rational agents whose strategy shares evolve under replicator dynamics. Henceforth, strategies earning above-average payoffs proliferate, while underperforming strategies decline. Recent theoretical work continues to extend the machinery of replicator dynamics to richer settings, including games in metric strategy spaces \citep{mendoza2024} and time-delayed payoffs \citep{varga2025}. 

In operations and policy contexts, EGT has proven particularly well suited to studying regulation because it captures how regulated firms gradually adapt to policy rather than jumping to a new equilibrium. \citet{long2019co} developed a three-party evolutionary game between government, consumers, and businesses to design waste management regulations. \citet{cheng2022egt} constructed a two-population evolutionary game of generator bidding in electricity markets. \citet{qian2026} analyzed the co-evolution of government environmental regulation and the adoption of green technologies by energy-intensive firms. 

In transportation and logistics specifically, \citet{zhangli2022} modeled the choice of mode of travel between car sharing and private cars using a logit dynamic evolutionary game. \citet{ahmad2023} provide a state-of-the-art review of EGT applications to urban road transport networks. \citet{zhang2025carbon} constructed a four-party evolutionary game between the government, ports and shipping companies to study the reduction of carbon emissions in maritime transport. \citet{norouzi2026using} analyzed trust, reputation, and pricing in gig-economy logistics platforms, modeling blockchain adoption as an evolutionary game. 

These studies collectively demonstrate the value of EGT for regulatory design. However, addressing the allocation of freight capacity or the regulation of quotas is a new topic in this field. Furthermore, studies that model a regulator using EGT assume the regulator's strategy space as an exogenously specified set of policy levers rather than a menu of uncertainty-robust, Pareto-efficient policies derived from an optimization model.

\subsection{Research Gaps and Positioning}

The above review exposes three gaps at the intersection of the addressed streams. First, the only game-theoretic treatment for regulating the minimum quota in freight transportation \citep{pazoki2024} is deterministic and static. In other words, it assumes perfect knowledge of demand, prices, and costs, treats enforcement as costless, and presumes instantaneous convergence to a perfectly rational equilibrium. In addition, although robust optimization has been widely applied to supply chain and transportation planning, connecting a robust regulator to a game-theoretic model of market response is a new approach to designing regulations, and more specifically, quota regulations. Finally, EGT methods have been applied to the regulation of environmental, energy and urban transport. However, its application to the regulation of rail freight quota, where a regulated population is modeled, is novel. Especially when strategy space is specified exogenously rather than derived from an optimization model. 

The framework developed in this paper addresses all three gaps jointly. The proposed components are presented in Sections~\ref{sec:problem} to~\ref{sec:egt}. Table~\ref{tab:lit-position} compares the current article with the most closely related works.

\begin{table}[H]
\centering
\caption{A comparison the current study with the most closely related studies.}
\label{tab:lit-position}
\renewcommand{\arraystretch}{1.3}
\resizebox{\textwidth}{!}{%
\begin{tabular}{lcccccc}
\toprule
\textbf{Study} & \textbf{Three players} & \textbf{Uncertainty} 
& \textbf{Multi-objective} & \textbf{Adaptive (EGT)} 
& \textbf{Quota regulation} & \textbf{Endogenous policy selection} \\
\midrule
\citet{amaruchkul2011}   &            & \checkmark &            &            &            &            \\
\citet{taherkhani2022}   &            & \checkmark &            &            &            &            \\
\citet{mirzaee2022}      & \checkmark & \checkmark &            &            &            &            \\
\citet{mirzaee2023}      &            & \checkmark & \checkmark &            &            &            \\
\citet{pazoki2024}       & \checkmark &            &            &            & \checkmark &            \\
\citet{zhang2025carbon}  &            &            &            & \checkmark &            &            \\
\citet{qian2026}         &            &            &            & \checkmark &            &            \\
\citet{norouzi2026using} &            &            &            & \checkmark &            &            \\
\textbf{The current paper}      & \checkmark & \checkmark & \checkmark & \checkmark & \checkmark & \checkmark \\
\bottomrule
\end{tabular}
}
\end{table}

\section{Problem Description and Assumptions}
\label{sec:problem}

\subsection{Problem Statement}
\label{subsec:problem-statement}

This paper studies the design of a minimum shipment quota regulation in a rail freight market characterized by limited transportation capacity, heterogeneous commodity margins, and parameter uncertainty. The setting is motivated by the North American freight system, in which a single rail carrier with limited annual capacity transports multiple commodities produced by distinct suppliers, and a public regulator can intervene to ensure that socially critical commodities are not excluded from the market by purely profit-driven capacity allocations.

We assume that three players participate in a strategic game in which each player pursues different goals. The regulator tries to achieve the highest possible social welfare through minimum shipment quotas imposed on commodities that otherwise would not be shipped under any circumstances. The implementation of such quotas is accompanied by administrative, supervisory, and political costs. Producers, who make individual commodities, fight for the capacity of the shipper horizontally by making bids regarding transport prices. Their goal is to achieve maximum net profits after deducting transport costs and the penalty for unshipped commodities. The shipper distributes its capacity among various commodities to maximize operating margin within the framework of quotas, if imposed, and the revenue constraint, which is called Maximum Revenue Entitlement (MRE) in Canada.

Two features distinguish this setting from the problem formulations addressed in the existing literature. First, key market primitives, most notably the volume of each commodity available for shipment, are uncertain at the time the regulator commits to a quota policy. Harvest yields, production levels, and export conditions are realized months after the regulatory framework has been established, so quotas must be determined here-and-now while shipment decisions are made adaptively once uncertainty is resolved. Second, convergence to a fully rational equilibrium is not instantaneous. Suppliers observe market outcomes over repeated periods and gradually adjust their bidding behavior, where the carrier experiments with alternative allocation rules, and the regulator revises its policy stance in response to observed market behavior. The question of interest is therefore not only which quota policy is optimal in a one-shot sense, but which policy remains evolutionarily stable as all three populations adapt.

Consequently, this paper investigates how minimum shipment quotas should be designed to maximize social welfare and remain robust to uncertainty in the volume of each commodity available for shipment, while explicitly accounting for the cost of enforcement. Furthermore, the article reviews, among the resulting set of robust, Pareto-efficient quota policies, the policy that emerges as evolutionarily stable when the suppliers and the carrier adapt their strategies over time.

This problem is addressed in two connected stages. In the first stage, an adaptive robust optimization (ARO) model characterizes the regulator's quota-setting problem under worst-case realizations of uncertain commodity availability, trading off social welfare against enforcement cost to generate a Pareto frontier of robust quota policies. In the second stage, this Pareto frontier defines the regulator's strategy space in an evolutionary game in which the regulator, the suppliers, and the carrier co-evolve under replicator dynamics. The mechanism through which the two stages are connected is presented in Section~\ref{sec:solution-framework}.

\subsection{Notation}\label{subsec:notation}

Table~\ref{tab:notation} summarizes the sets, parameters, decision variables, and evolutionary-game notation used throughout the paper. Commodities are indexed by $i \in \mathcal{I} = \{1, \ldots, n\}$. Two derived quantities recur throughout the analysis. The first is the carrier's per-unit transport margin on flow $i$, $p_i - c_i$. The second is the regulator's cost of enforcing quota $L_i$,
\begin{equation}
K_i(L_i) = \rho_i \mathbf{1}\{L_i > 0\} + \lambda_i L_i + \tfrac{1}{2}\eta_i L_i^2,
\label{eq:K-def}
\end{equation}
whose three terms are explained in Assumption~6. Symbols used only in specific reformulations are defined locally where they appear.

\begin{table}[!htbp]
\centering
\footnotesize
\caption{Notation used throughout the integrated framework.}
\label{tab:notation}
\renewcommand{\arraystretch}{1.05}
\setlength{\tabcolsep}{4pt}
\begin{tabularx}{\textwidth}{>{\centering\arraybackslash}p{0.155\textwidth} X}
\toprule
\textbf{Symbol} & \textbf{Description} \\
\midrule
\multicolumn{2}{l}{\textit{Sets and indices}} \\
$\mathcal{I}$ & Set of commodities, indexed by $i$. \\
$\mathcal{J}_i$ & Set of admissible bidding strategies for supplier $i$, indexed by $j$. \\
$\mathcal{A}$ & Set of admissible carrier allocation policies, indexed by $a$. \\
$\mathcal{R}$ & Set of Pareto-efficient regulator quota policies, indexed by $r$. \\
$\mathcal{U}$ & Set of admissible realizations of the uncertain parameters. \\
$k$ & Index of iterations in the $\epsilon$-constraint procedure. \\
$t$ & Evolutionary time. \\
\addlinespace
\multicolumn{2}{l}{\textit{Parameters}} \\
$D_i$ & Available supply of commodity $i$. \\
$T$ & Total transportation capacity available to the rail carrier. \\
$S_i$ & Exogenous benchmark market price of commodity $i$. \\
$A_i$, $\kappa_i$ & Scale parameters of the isoelastic demand and inverse-demand functions. \\
$\varepsilon_i$ & Absolute price elasticity of demand for commodity $i$. \\
$c_i$ & Carrier's per-unit operating cost for transporting commodity $i$. \\
$g_i$ & Per-unit penalty on the quantity of commodity $i$ that is not shipped. \\
$\beta_i$ & Social value assigned by the regulator to one transported unit of commodity $i$. \\
$m_i$ & Carrier's transport margin at the maximum rational bid, $S_i + g_i - c_i$. \\
$\mathrm{MRE}_i$ & Maximum revenue entitlement associated with commodity $i$. \\
$\rho_i,\ \lambda_i,\ \eta_i$ & Fixed, linear, and quadratic enforcement-cost coefficients for commodity $i$. \\
$\varphi$ & Worst-case availability factor, $\underline{D}_i = \varphi D_i$. \\
$K^{\max}$ & Enforcement budget required to protect every quota-dependent commodity in full. \\
$\epsilon^k$ & Enforcement-cost budget in iteration $k$ of the $\epsilon$-constraint procedure. \\
$\mu$ & Mutation rate in the replicator dynamics. \\
$T_{\max}$ & Evolutionary horizon. \\
\addlinespace
\multicolumn{2}{l}{\textit{Uncertainty}} \\
$\boldsymbol{\xi}$ & Vector representing a realization of the uncertain parameters. \\
$\widetilde{\theta}$; $\underline{\theta}, \overline{\theta}$ & Generic uncertain parameter in $\boldsymbol{\xi}$, and its lower and upper bounds. \\
$\underline{D}_i$ & Worst-case available volume of commodity $i$ under the uncertainty set. \\
$u_i(\boldsymbol{\xi})$ & Adaptive shipment quantity of commodity $i$ under realization $\boldsymbol{\xi}$. \\
$\Phi(\boldsymbol{L},\boldsymbol{\xi})$ & Optimal value of the second-stage shipment problem. \\
\addlinespace
\multicolumn{2}{l}{\textit{Decision variables and market outcomes}} \\
$L_i$, $\boldsymbol{L}$ & Minimum shipment quota on commodity $i$, and the quota vector. \\
$L_i^{(r)}$, $\boldsymbol{L}^{(r)}$ & Quota for commodity $i$, and full quota vector, under regulator policy $r$. \\
$u_i$, $\boldsymbol{u}$ & Quantity of commodity $i$ transported by the carrier, and the shipment vector. \\
$p_{ij}$ & Transport-price bid associated with strategy $j$ of supplier $i$. \\
$\bar{p}_i$ & Mean-field bid of supplier $i$. \\
$Q_i$ & Quantity of commodity $i$ supplied to its market; $Q_i = u_i$. \\
$K_i(L_i)$ & Enforcement cost of imposing quota $L_i$ on commodity $i$. \\
\addlinespace
\multicolumn{2}{l}{\textit{Evolutionary-game notation}} \\
$x_{ij}(t)$ & Share of supplier-$i$ agents adopting bidding strategy $j$ at time $t$. \\
$y_a(t)$ & Share of carrier agents adopting allocation policy $a$ at time $t$. \\
$z_r(t)$ & Share of regulator agents adopting quota policy $r$ at time $t$. \\
$\pi_{ij}$, $\overline{\pi}_i$ & Expected payoff of supplier-$i$ strategy $j$, and supplier-$i$ population average. \\
$\Pi_a^{\mathrm{car}}$, $\overline{\Pi}^{\mathrm{car}}$ & Expected carrier payoff from policy $a$, and carrier population average. \\
$U_r^{\mathrm{reg}}$, $\overline{U}^{\mathrm{reg}}$ & Expected regulator utility from policy $r$, and regulator population average. \\
\bottomrule
\end{tabularx}
\end{table}

\subsection{Assumptions}
\label{subsec:assumptions}

The model is developed under the following assumptions.

\begin{enumerate}
    \item \textbf{Tri-level Stackelberg timing.} Decisions unfold hierarchically. The regulator moves first by setting minimum shipment     quotas. The suppliers move second by simultaneously posting transport price bids, and the carrier moves last by allocating capacity. Each player anticipates the reactions of the next players in the hierarchy.

    \item \textbf{Isoelastic demand.} Each commodity $i$ faces a constant price-elasticity demand curve $Q_i(S_i) = A_i S_i^{-\varepsilon_i}$, with inverse demand $S_i(Q_i) = \kappa_i Q_i^{-1/\varepsilon_i}$, where $A_i, \kappa_i, \varepsilon_i > 0$. The demand intercept $A_i$ is calibrated to the observed price--quantity pairs. The isoelastic form governs the evolutionary stage, while the robust stage uses a linear approximation for tractability.

    \item \textbf{Uncertain commodity availability;} The volume of each commodity available for shipment is uncertain at the time quota are set.  Uncertainty is represented by an interval (box) uncertainty set around market prices and customer demand, and the regulator adopts a worst-case (ambiguity-averse) posture consistent with adaptive robust optimization.

    \item \textbf{Bounded continuous knapsack carrier.} The carrier's allocation problem is a bounded continuous knapsack with hard lower bounds. The mandatory quota $L_i$ is first honored, after which the residual capacity is allocated according to the carrier's prevailing policy, up to the available volume of each commodity. A profit-oriented carrier does not voluntarily ship commodities whose transport margin is negative beyond the mandated quota.

    \item \textbf{Revenue cap (MRE).} A maximum revenue entitlement cap limits the carrier's transport revenue on each flow, $u_i(p_i - c_i) \le \mathrm{MRE}_i$, which in turn bounds the feasible supplier bids from above.

    \item \textbf{Separable, convex enforcement cost.} The regulator's cost of enforcing quota $L_i$ is separable between commodities and takes the form of Eq.~\eqref{eq:K-def}, comprising a fixed activation cost, a linear per-unit oversight burden, and a convex rising marginal burden.

    \item \textbf{Finite strategy discretization;} for evolutionary analysis, each player's strategy space is discretized into a finite set. Supplier $i$ selects from a finite price grid $\{p_{ij}\}_{j \in J_i}$, the carrier selects from a finite set $\mathcal{A}$ of allocation policies, and the regulator selects from a finite menu $\mathcal{R}$ of Pareto-efficient quota vectors generated by the ARO model.

    \item \textbf{Replicator dynamics with mutation.} Strategy shares evolve according to discrete-time replicator dynamics. Strategies earning above-average payoffs increase in population share, while strategies earning below-average payoffs decline. A small mutation rate maintains the diversity of the strategy and prevents premature extinction.

    \item \textbf{Single aggregate rail capacity.} The rail network is represented by a single annual capacity constraint $\sum_{i \in \mathcal{I}} u_i \le T$, abstracting from route-level and seasonal congestion.

    \item \textbf{Single market per commodity.} Each commodity clears in a single destination market, so that the shipped quantity coincides with the marketed quantity ($Q_i = u_i$), and the market price $S_i$ is treated as exogenous to any individual supplier's bidding decision.
\end{enumerate}

\section{Solution Framework}\label{sec:solution-framework}

This section specifies how robust optimization and evolutionary stages are combined. The mechanism connecting the two models is presented first, followed by the complete solution procedure and its formal statement as an algorithm.

\subsection{The Integration Mechanism} \label{subsec:integration}

The ARO and EGT models do not operate separately. The results of the ARO model define the strategy space for the players participating in the EGT model. These two models provide answers to different but related questions, and their combination is a contribution made from a methodological perspective in this research. The ARO model identifies the quota regulation the regulator may impose considering the uncertainty of commodity supplies and the conflict between welfare maximization and cost minimization. As bigger quotas lead to higher welfare but also higher costs, the optimal quota vector does not exist. Instead, several Pareto-efficient vectors are obtained. The EGT model identifies which of the Pareto-efficient quota vectors will be dominant when strategic and adaptive interactions occur between regulator, suppliers, and carrier.

The connection between the two operates through the Pareto frontier generated by the ARO stage. Let
\begin{equation}
\mathcal{R} = \left\{ L^{(1)}, L^{(2)}, \ldots, L^{(R)} \right\}
\label{eq:pareto-set}
\end{equation}
denote the set of Pareto-efficient and uncertainty-robust quota vectors obtained by solving the ARO model over a grid of enforcement-cost budgets (the $\varepsilon$-constraint method is detailed in Section~\ref{sec:aro}). Each element $L^{(r)} = (L^{(r)}_1, \ldots, L^{(r)}_n)$ is a complete quota policy specifying a minimum shipment level for every commodity, and each is Pareto non-dominated in the sense that no other feasible policy simultaneously delivers higher worst-case social welfare and lower enforcement cost.

In the evolutionary stage, $\mathcal{R}$ is adopted directly as the regulator's strategy space. The regulator population does not evolve over the continuum of all conceivable quota vectors in $[0, D_1] \times \cdots \times [0, D_n]$, but only over the $R$ Pareto-efficient policies contained in $\mathcal{R}$. The population share $z_r$ denotes the fraction of the regulator population that adopts the policy $L^{(r)}$, and the vector $z = (z_1, \ldots, z_R)$ evolves under replicator dynamics together with the strategy distributions of the suppliers and the carrier.

This construction delivers two key findings of the framework. First, any quota policy analyzed in the evolutionary stage is uncertainty-robust and Pareto efficient by definition. Therefore, the evolutionary dynamics does not assign any probabilities to dominated and/or fragile policies. Second, the strategy space of the regulator becomes finite (cardinality $R$), instead of being continuous, making the computation of the coupled replicator system possible and delivering a full range of tradeoffs between efficiency and enforcement at the same time. In such a manner, the ARO stage ensures the quality and feasibility of regulator options, while the EGT stage governs their choice due to adaptation to strategic environment.

This functional segregation ensures the higher efficacy of the coupled framework compared to its components used separately. In particular, a robust optimization model, when used alone, provides a Pareto frontier but lacks a systematic approach to choosing one option as it usually requires a separate decision maker or arbitrary objective weighing. An evolutionary game, when used separately, has a built-in mechanism for choosing a stable solution but does not guaranty that the underlying strategies are both Pareto efficient and uncertainty-robust. Coupling of the two models ensures that evolutionary dynamics will take place among a properly defined set of Pareto efficient and uncertainty-robust quota policies and the selected policy is uncertainty-robust, evolutionarily stable, and Pareto efficient at the same time.

The connection in question is unidirectional in construction but bidirectional in effect. The Pareto frontier is calculated before entering the evolutionary stage and is not re-evaluated during the dynamics, so the output of the ARO stage flows in the EGT model. However, the prevailing policy is chosen by the interaction of all three players, which means that the bidding behavior of suppliers and allocation behavior of carrier feed back into the relative fitness of quota policies in question. Consequently, the equilibrium policy is not the optimal policy in terms of welfare but is evolutionarily stable taking into account other players' strategic reactions.

\subsection{Solution Algorithm}
\label{subsec:solution-algorithm}

The complete solving mechanism consists of two consecutive stages: a robust optimization stage that creates the strategy space of the regulator, and an evolutionary stage that defines the existing equilibrium. Both stages are described by Algorithm~\ref{alg:aro-egt}.

\paragraph{Phase 1: Construction of the robust Pareto frontier.}
The ARO model is solved repeatedly on a grid of $R$ evenly spaced enforcement budgets spanning $[0, K^{\max}]$, where $K^{\max}$ is the maximum feasible enforcement expenditure subject to capacity and availability constraints. Each instance is solved under the worst-case realization of commodity availability, $\underline{D}_i = \varphi D_i$, and is warm-started from the solution of the preceding budget level, which promotes numerical stability and continuity of the frontier. The resulting set $\mathcal{R} = \{L^{(1)}, \ldots, L^{(R)}\}$ constitutes the regulator's strategy menu.

\paragraph{Phase 2: Evolutionary selection.}
The three populations co-evolve over their respective strategy spaces: the suppliers over their price-bid grids, the carrier over its allocation policies $\mathcal{A}$, and the regulator over the Pareto menu $\mathcal{R}$. All are initialized with uniform shares, reflecting the absence of prior information on which strategies will prevail. At each step, payoffs are evaluated as mean-field expectations over the current strategy distributions of the other two players, so that each player responds best to the average behavior of its counterparts rather than to any single realization. Shares are then updated by the discrete-time replicator rule with mutation. The process terminates once the distributions stabilize or the horizon $T_{\max}$ is reached, at which point the dominant strategy of each population identifies the equilibrium quota policy $L^{\ast}$, supplier bids $p^{\ast}$, and carrier policy $a^{\ast}$, from which the shipments $u^{\ast}$, payoffs, and social welfare achieved follow.

\begin{algorithm}[H]
\caption{Integrated ARO--EGT Solution Procedure}
\label{alg:aro-egt}
\begin{algorithmic}[1]
\Require Commodity data $\{D_i, S_i, c_i, g_i, \varepsilon_i, \beta_i, \rho_i, \lambda_i, \eta_i\}_{i\in\mathcal{I}}$; capacity $T$; worst-case factor $\varphi$; frontier size $R$; carrier policies $\mathcal{A}$; price grids $\{J_i\}$; mutation rate $\mu$; horizon $T_{\max}$
\Ensure Equilibrium quota policy $L^{\ast}$, bids $p^{\ast}$, allocation policy $a^{\ast}$, shipments $u^{\ast}$

\Statex \textbf{Phase 1: Robust Pareto frontier}
\State $\underline{D}_i \gets \varphi\, D_i$ for all $i\in\mathcal{I}$ 
\Comment{worst-case availability}
\State $K^{\max} \gets$ maximum feasible enforcement cost s.t. capacity and availability
\State Partition $[0, K^{\max}]$ into budgets 
$B_1 < \cdots < B_R$
\State $L^{(0)} \gets \varnothing$ \Comment{warm-start seed}
\For{$r = 1$ to $R$}
    \State $L^{(r)} \gets$ solve ARO $\varepsilon$-constraint problem at budget $B_r$, warm-started from $L^{(r-1)}$
\EndFor
\State $\mathcal{R} \gets \{L^{(1)}, \ldots, L^{(R)}\}$ 
\Comment{regulator strategy menu}

\Statex \textbf{Phase 2: Evolutionary co-evolution}
\State Initialize $x_{ij} \gets 1/|J_i|$, $y_a \gets 1/|\mathcal{A}|$, $z_r \gets 1/R$ \Comment{uniform shares}
\For{$t = 1$ to $T_{\max}$}
    \State $\bar{p}_i \gets \sum_{j\in J_i} x_{ij}\, p_{ij}$ for all $i$ 
    \Comment{mean-field bids}
    \For{each regulator policy $r$ and carrier policy $a$}
        \State $u \gets \textsc{Allocate}(L^{(r)}, \bar{p}, a)$
        \State accumulate carrier payoff $\Pi^{\mathrm{car}}_a$ and regulator payoff $U^{\mathrm{reg}}_r$
        \For{each supplier $i$ and bid $j \in J_i$}
            \State $u^{(j)} \gets \textsc{Allocate}(L^{(r)}, (p_{ij}, 
            \bar{p}_{-i}), a)$ \Comment{$i$ deviates, others at mean}
            \State $\pi_{ij} \gets \pi_{ij} + y_a z_r\, 
            \pi_i(u^{(j)}_i, p_{ij})$
        \EndFor
    \EndFor
    \State $x_{ij} \gets \textsc{ReplicatorUpdate}(x_{ij}, \pi_{ij}, \mu)$ 
    for all $i,j$
    \State $y_a \gets \textsc{ReplicatorUpdate}(y_a, \Pi^{\mathrm{car}}_a, 
    \mu)$ for all $a$
    \State $z_r \gets \textsc{ReplicatorUpdate}(z_r, U^{\mathrm{reg}}_r, \mu)$ 
    for all $r$
\EndFor

\Statex \textbf{Output}
\State $a^{\ast}\!\gets\!\arg\max_a y_a$;\ \ 
$r^{\ast}\!\gets\!\arg\max_r z_r$;\ \ 
$p^{\ast}_i\!\gets\! p_{i j^{\ast}}$ with 
$j^{\ast}\!=\!\arg\max_j x_{ij}$
\State $L^{\ast} \gets L^{(r^{\ast})}$;\ \ 
$u^{\ast} \gets \textsc{Allocate}(L^{\ast}, p^{\ast}, a^{\ast})$
\State \Return $L^{\ast}, p^{\ast}, a^{\ast}, u^{\ast}$
\end{algorithmic}
\end{algorithm}

The subroutine $\textsc{Allocate}(L, p, a)$ implements the carrier's bounded continuous knapsack, detailed in Section~\ref{subsec:allocation}, and $\textsc{ReplicatorUpdate}$ applies the discrete-time replicator rule with mutation of Section~\ref{subsec:replicator}.

\begin{figure}[H]
\centering
\resizebox{\textwidth}{!}{%
\begin{tikzpicture}[
    node distance = 1.0cm and 1.0cm,
    box/.style   = {rectangle, rounded corners, draw=black, thick,
                    align=center, minimum height=1.15cm, minimum width=2.9cm,
                    fill=gray!8, font=\footnotesize},
    datastore/.style = {rectangle, draw=black, thick, align=center,
                    minimum height=1.15cm, minimum width=2.9cm, fill=blue!7,
                    font=\footnotesize},
    phase/.style = {rectangle, draw=black!55, thick, dashed,
                    rounded corners, inner sep=0.4cm},
    arrow/.style = {-{Latex[length=2.4mm]}, thick},
    lbl/.style   = {font=\footnotesize\itshape, text=black!65, align=center}
]

\node[datastore] (data) {Commodity data \& \\ uncertainty set $\mathcal{U}$};
\node[box, right=of data] (aro) {ARO model \\ 
$\max_{\boldsymbol L}\min_{\boldsymbol\xi\in\mathcal U}$ \\ (two-stage)};
\node[box, right=of aro] (rc) {Robust counterpart \\ worst-case reduction};
\node[box, right=of rc] (pareto) {Pareto frontier \\ 
$\mathcal{R}=\{L^{(1)},\dots,L^{(R)}\}$};

\node[box, below=2.2cm of pareto] (reg) {Regulator \\ evolves over 
$\mathcal{R}$};
\node[box, left=of reg] (sup) {Suppliers \\ evolve over bids};
\node[box, left=of sup] (car) {Carrier \\ evolves over $\mathcal{A}$};
\node[datastore, left=of car] (eq) {Equilibrium \\ 
$L^{\ast},\,p^{\ast},\,a^{\ast},\,u^{\ast}$};

\draw[arrow] (data) -- (aro);
\draw[arrow] (aro)  -- (rc);
\draw[arrow] (rc) -- (pareto);
\node[lbl, inner sep=1.5pt, above=0.30cm] (epslbl) 
    at ($(rc)!0.5!(pareto)+(0,0.575)$) {$\varepsilon$-constraint \\ sweep};
\draw[black!45, thin] (epslbl.south) -- ($(rc)!0.5!(pareto)$);

\draw[arrow] (pareto) -- node[lbl, right] {strategy \\ menu} (reg);

\draw[arrow, <->] (reg) -- (sup);
\draw[arrow, <->] (sup) -- (car);

\draw[arrow] (car) -- (eq);

\node[below=0.9cm of sup, inner sep=0pt, minimum size=0pt] (dip) {};
\draw[arrow, <->] (reg.south) to[out=-90,in=0] (dip.center) 
                              to[out=180,in=-90] (car.south);

\begin{scope}[on background layer]
    \node[phase, fit=(data)(aro)(rc)(pareto)(epslbl)] (p1frame) {};
    \node[phase, fit=(eq)(car)(sup)(reg)(dip)] (p2frame) {};
\end{scope}

\node[lbl, above=0.12cm of p1frame] {Phase 1: Robust optimization (ARO)};
\node[lbl, below=0.12cm of p2frame] {Phase 2: Evolutionary selection (EGT)};

\end{tikzpicture}
}
\caption{The integrated ARO--EGT solution framework. \\In Phase~1, the two-stage adaptive robust optimization model is reduced to its deterministic robust counterpart and solved repeatedly under the $\varepsilon$-constraint method, generating a Pareto frontier of robust, Pareto-efficient quota policies. \\In Phase~2, this frontier serves as the regulator's strategy menu, and the regulator, suppliers, and carrier co-evolve under replicator dynamics (double-headed arrows) until the equilibrium $(L^{\ast}, p^{\ast}, a^{\ast}, u^{\ast})$ is reached.}
\label{fig:framework}
\end{figure}
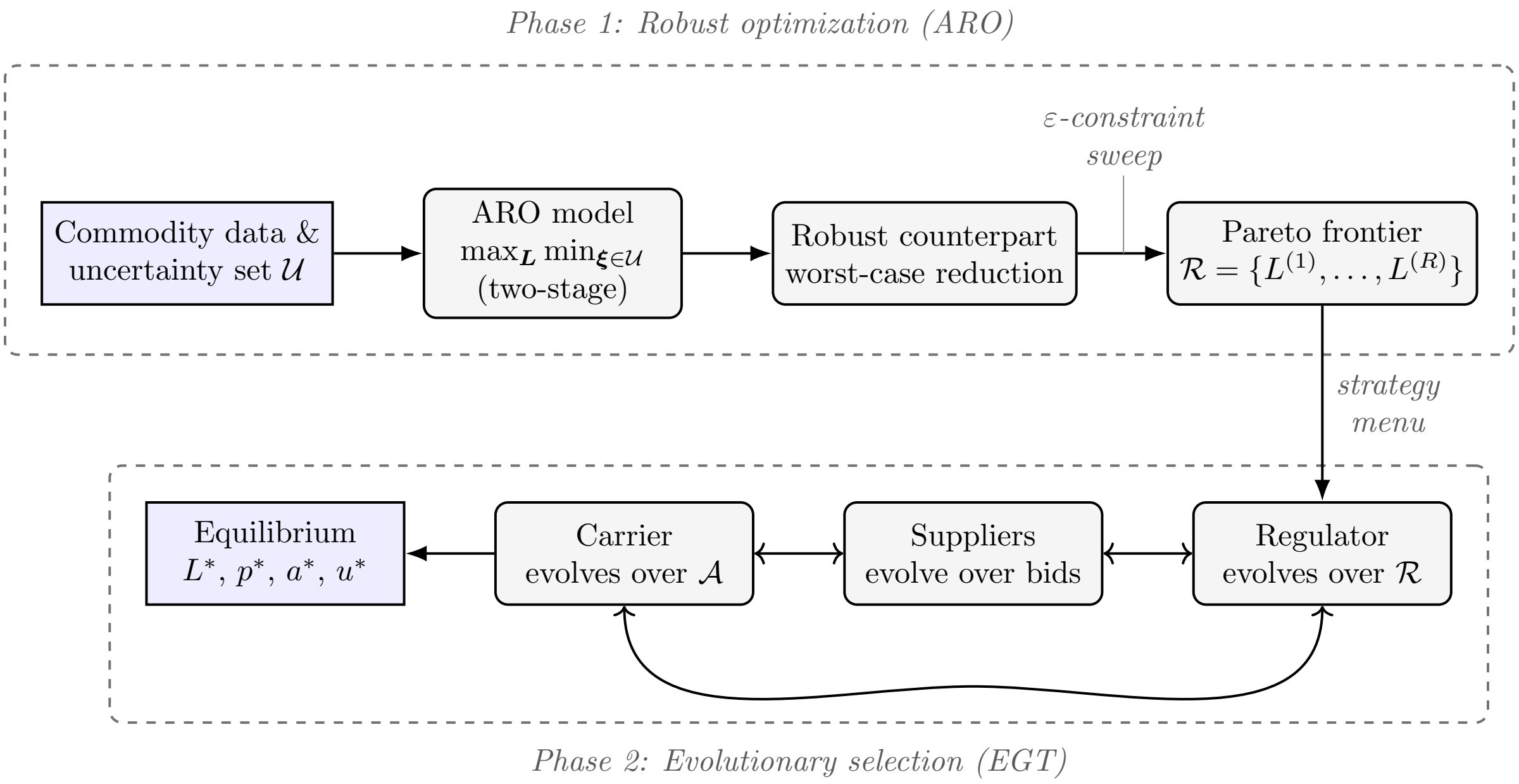

Figure~\ref{fig:framework} summarizes the integrated solution framework. The robust optimization stage (upper block) transforms the commodity data and uncertainty set into a Pareto frontier of quota policies. This frontier is passed to the evolutionary stage (lower block), where it serves as the regulator's strategy menu, and the three populations co-evolve to the prevailing equilibrium.

\section{The Adaptive Robust Optimization Model (Regulator's Problem)}
\label{sec:aro}

The regulator must determine minimum shipment quotas before uncertain market conditions are observed. Once uncertainty is realized, shipment quantities are adjusted subject to transportation-capacity limitations, commodity availability, and mandated quotas. The quota-design problem therefore has a natural two-stage structure, which an Adaptive Robust Optimization (ARO) framework is designed to capture, in the sense that quotas are here-and-now decisions, while shipment quantities are adaptive decisions taken after uncertainty is resolved. The model identifies quota policies that remain effective under adverse market conditions and generates a set of robust Pareto-efficient solutions, which form the regulator's strategy space in the evolutionary game of Section~\ref{sec:egt}.

\subsection{Model Setting and Uncertainty}
\label{subsec:aro-setting}

\paragraph{Decision variables and feasibility.}
The regulator determines a quota vector $\boldsymbol L=(L_i)_{i\in\mathcal I}$, where $L_i$ is the minimum shipment quota imposed on commodity $i$; the adaptive shipment quantities are $\boldsymbol u=(u_i)_{i\in\mathcal I}$. Both are bounded by the volume available for shipment, and total shipments are limited by rail capacity:
\begin{equation}
0 \le L_i \le u_i \le D_i(S_i), \quad \forall i\in\mathcal I,
\qquad
\sum_{i\in\mathcal I} u_i \le T.
\label{eq:aro_feasibility}
\end{equation}
The chain $L_i \le u_i$ is the quota-enforcement requirement. Hence, shipments must meet or exceed the minimum required. Note that \eqref{eq:aro_feasibility} implies $\sum_{i} L_i \le T$, so a separate feasibility restriction is not required on the quota vector.

\paragraph{Uncertainty set.}
Probability distributions are not assumed to be known. Instead, the market price of each commodity is taken to lie in a bounded interval,
\begin{equation}
\underline S_i \le \widetilde S_i \le \overline S_i,
\qquad \forall i\in\mathcal I,
\label{eq:price_bounds}
\end{equation}
so that the uncertainty set is
\begin{equation}
\mathcal U = \Big\{ \boldsymbol\xi = (\widetilde S_i)_{i\in\mathcal I} :
\underline S_i \le \widetilde S_i \le \overline S_i,
\; \forall i\in\mathcal I \Big\}.
\label{eq:aro_uncertainty}
\end{equation}

\paragraph{Demand representation.}
Because the ARO stage serves to identify robust quota policies rather than to represent strategic market interactions, demand is considered to be linear in price,
\begin{equation}
D_i(S_i) = a_i - b_i S_i, \qquad a_i>0, \; b_i>0,
\label{eq:linear_demand}
\end{equation}
so that higher realized prices reduce the volume available for transportation. The isoelastic specification of Assumption~2 governs the evolutionary stage; the linear form is adopted here to keep the robust counterpart tractable.

\subsection{The Regulator's Robust Problem}
\label{subsec:aro-problem}

\paragraph{Carrier behavior.}
Although the regulator sets the quota, it does not supervise the process of allocating the residual capacity. The profit-oriented carrier will move a commodity in excess of its quota, provided that the process results in a positive transportation margin (Assumption 4). The fact that the ARO stage comes before the evolutionary game implies that the bids have not yet been made. Therefore, the carrier’s margin is calculated based on the maximum bid that a supplier can rationally make:
\begin{equation}
m_i = S_i + g_i - c_i.
\label{eq:margin}
\end{equation}
Commodities with $m_i < 0$ are never carried voluntarily:
\begin{equation}
u_i \le L_i, \qquad \forall i \in \mathcal I \ \text{with} \ m_i < 0.
\label{eq:no_voluntary}
\end{equation}
Together with $u_i \ge L_i$ in \eqref{eq:aro_feasibility}, condition
\eqref{eq:no_voluntary} forces $u_i = L_i$ for this commodity, reaching the market in the exact mandated volume and no more. This makes the quota regulation binding. In the evolutionary stage of Section~\ref{sec:egt}, the carrier responds instead to the prevailing bids $p_i$ as they emerge.

\paragraph{Enforcement cost.}
The cost of enforcing quota $L_i$ is $K_i(L_i)$ as defined in Eq.~\eqref{eq:K-def}, totaling
\begin{equation}
K(\boldsymbol L) = \sum_{i\in\mathcal I} K_i(L_i).
\label{eq:aro_total_cost}
\end{equation}
The fixed activation component $\rho_i$ is set to zero throughout, so that $K_i(L_i) = \lambda_i L_i + \tfrac{1}{2}\eta_i L_i^2$ and $K(\boldsymbol L)$ is convex and quadratic.

\paragraph{Robust objective.}
The two-stage structure is made explicit through the second-stage value function. For a fixed quota vector $\boldsymbol L$ and a realized $\boldsymbol\xi$, the shipment decision solves:
\begin{equation}
\Phi(\boldsymbol L, \boldsymbol\xi) = \max_{\boldsymbol u} \;
\sum_{i\in\mathcal I} \beta_i u_i
\quad \text{s.t.} \quad
\sum_{i\in\mathcal I} u_i \le T, \;\;
L_i \le u_i \le D_i(\boldsymbol\xi), \;\;
u_i \le L_i \ \text{if} \ m_i < 0.
\label{eq:recourse_value}
\end{equation}
The regulator, who is ambiguity-averse and commits to $\boldsymbol L$ before $\boldsymbol\xi$ solves:
\begin{equation}
\max_{\boldsymbol L} \; \min_{\boldsymbol\xi\in\mathcal U} \;
\Phi(\boldsymbol L, \boldsymbol\xi) - K(\boldsymbol L) \quad \text{s.t.} \quad 0 \le L_i \le \min_{\boldsymbol\xi\in\mathcal U} D_i(\boldsymbol\xi),
\label{eq:aro_original}
\end{equation}
which is required for the second stage to remain feasible under every realization.

\subsection{The Deterministic Robust Counterpart}
\label{subsec:aro-counterpart}

The adjustable problem \eqref{eq:aro_original} admits an exact deterministic reduction. Availability enters \eqref{eq:recourse_value} only through the upper bounds $u_i \le D_i(\boldsymbol\xi)$, so raising any $D_i$ enlarges the second-stage feasible set and cannot decrease $\Phi$. Therefore, the value function is non-decreasing in $\boldsymbol D(\boldsymbol\xi)$ componentwise. Since $D_i(S_i)$ is decreasing in price, the least favorable availability occurs at the largest admissible price, and the inner minimum is attained at the componentwise lower bound:
\begin{equation}
\underline D_i = a_i - b_i \overline S_i, \qquad \forall i\in\mathcal I,
\label{eq:worst_case_demand}
\end{equation}
hence, $\min_{\boldsymbol\xi\in\mathcal U} \Phi(\boldsymbol L,\boldsymbol\xi) = \Phi(\boldsymbol L, \underline{\boldsymbol D})$. Because the worst case is associated with a specific realization rather than a probability distribution over several realizations, the recourse decision has to be optimal only for that particular realization, resulting in a simplification of the two-stage model into a one-stage problem. Hence, adaptivity must be seen as a characteristic of model timing but not as another source of flexibility in this situation. Adaptivity determines what variables are bound by which realization. However, when combined with a box uncertainty set and a monotonic recourse value, it results in the same strategy as minimization of the worst-case objective. 

The price bounds are calibrated so that this worst case corresponds to a uniform proportional shortfall,
\begin{equation}
\underline D_i = \varphi D_i, \qquad \varphi \in (0,1],
\label{eq:phi_def}
\end{equation}
which is the form used throughout the remainder of the paper. The baseline value $\varphi = 0.9$ is discussed in Section~\ref{subsec:calibration}. Substituting \eqref{eq:phi_def} into \eqref{eq:aro_original} eliminates the inner minimization and yields a deterministic single-level program:
\begin{equation}
\begin{aligned}
\max_{\boldsymbol L, \boldsymbol u} \quad
& \sum_{i\in\mathcal I} \beta_i u_i - K(\boldsymbol L) \\
\textnormal{s.t.} \quad
& \sum_{i\in\mathcal I} u_i \le T, \\
& L_i \le u_i, \qquad \forall i\in\mathcal I, \\
& u_i \le L_i, \qquad \forall i\in\mathcal I \ \text{with} \ m_i < 0, \\
& 0 \le L_i \le \underline D_i, \quad
  0 \le u_i \le \underline D_i, \qquad \forall i\in\mathcal I.
\end{aligned}
\label{eq:robust_counterpart}
\end{equation}
All constraints are linear and the objective is concave, so \eqref{eq:robust_counterpart} is a convex program that can be solved directly with a standard optimization package.

\subsection{Bi-objective Formulation and the Pareto Frontier}\label{subsec:aro-pareto}

The regulator balances two competing objectives, i.e., maximizing worst-case social welfare:
\begin{equation}
f_1(\boldsymbol L, \boldsymbol u) = \sum_{i\in\mathcal I} \beta_i u_i,
\label{eq:f1}
\end{equation}
and minimizing enforcement cost:
\begin{equation}
f_2(\boldsymbol L) = K(\boldsymbol L).
\label{eq:f2}
\end{equation}
By applying the $\epsilon$-constraint method \citep{berube2009exact}, welfare is maximized subject to a bound on enforcement expenditure:
\begin{equation}
\max_{\boldsymbol L, \boldsymbol u} \; f_1(\boldsymbol L, \boldsymbol u)
\quad \text{subject to} \quad
K(\boldsymbol L) \le \epsilon^k,
\label{eq:eps_constraint}
\end{equation}
together with constraints \eqref{eq:robust_counterpart}. The quadratic form of $K$ makes \eqref{eq:eps_constraint} a quadratically
constrained program that remains convex.

The upper end of the budget grid is set by the largest enforcement expenditure that is feasible under capacity and availability,
\begin{equation}
K^{\max} = \max \Big\{ K(\boldsymbol L) \; : \;
0 \le L_i \le \underline D_i, \;
\sum_{i\in\mathcal I} L_i \le T \Big\},
\label{eq:kmax}
\end{equation}
and the interval $[0, K^{\max}]$ is discretized into $R$ evenly spaced levels. Solving \eqref{eq:eps_constraint} at each level, warm-started from the preceding solution, generates the set of Pareto-efficient quota vectors:
\begin{equation}
\mathcal R = \Big\{ \boldsymbol L^{(1)}, \boldsymbol L^{(2)}, \dots,
\boldsymbol L^{(R)} \Big\}.
\label{eq:pareto_frontier}
\end{equation}
Each element of $\mathcal R$ is a robust and Pareto-efficient quota policy that balances transportation welfare against enforcement cost. This set is the regulator's strategy space in the evolutionary game, as established in Section~\ref{subsec:integration}.

\section{The Evolutionary Game-Theoretic Model}\label{sec:egt}

This section presents the EGT model that governs the second phase of the solution framework. The players and their strategy spaces are introduced first, followed by the demand structure, the shared feasibility constraints, the payoff functions, the carrier's allocation procedure, and finally the replicator dynamics through which the three populations co-evolve.

\subsection{Players, Strategy Spaces, and Timing}\label{subsec:players}

The model involves three players, the regulator, the suppliers, and the carrier, whose interaction follows the tri-level Stackelberg timing introduced in Assumption~1. The regulator moves first by committing to a minimum shipment quota vector, the suppliers move second by simultaneously posting transport price bids, and the carrier moves last by allocating its capacity in response to the observed quotas and bids. In the evolutionary setting, each player is represented not as a single optimizing agent but as a population of boundedly rational agents, and the distribution of strategies within each population evolves over time.

\paragraph{Regulator.}
The regulator selects a quota vector from the finite menu $\mathcal{R} = \{L^{(1)}, \ldots, L^{(R)}\}$ of robust, Pareto-efficient policies produced by the ARO stage, as established in Section~\ref{subsec:integration}. The population share $z_r \in [0,1]$ denotes the fraction of the regulator population adopting policy $L^{(r)}$, with $\sum_{r=1}^{R} z_r = 1$. Restricting the regulator's strategy space to $\mathcal{R}$ ensures that only efficient and uncertainty-robust policies participate in the evolutionary dynamics.

\paragraph{Suppliers.}
Each supplier $i \in \mathcal{I}$ selects a transport price bid from a finite grid 
\begin{equation}
J_i = \left\{ p_{i1}, p_{i2}, \ldots, p_{i,|J_i|} \right\} \subseteq \left[ 0,\; S_i + g_i \right],
\label{eq:price-grid}
\end{equation}
where the upper bound $S_i + g_i$ is the maximum bid a rational supplier would ever receive, since beyond this level the supplier prefers to store its product or to entirely abandon shipment \citep{pazoki2024}. The population share $x_{ij} \in [0,1]$ denotes the fraction of supplier $i$ agents who post bid $p_{ij}$, with $\sum_{j \in J_i} x_{ij} = 1$. The grid is taken to be evenly spaced over $[0, S_i + g_i]$, with $|J_i| = 10$ points in the case study.

\paragraph{Carrier.}
The carrier selects an allocation policy from a finite set $\mathcal{A} = \{a_0, a_1, \ldots\}$, where each policy is a rule to distribute residual capacity among commodities once the mandated quota have been met. The set $\mathcal{A}$ spans behaviorally distinct rules, ranging from purely profit-maximizing allocation to socially oriented and balanced allocation; the specific policies employed in the case study are defined in Section~\ref{subsec:allocation}. The population share $y_a \in [0,1]$ denotes the fraction of the carrier population using policy $a$, with $\sum_{a \in \mathcal{A}} y_a = 1$.

Collectively, the state of the evolutionary system at any time is described by the strategy distributions $x = (x_{ij})$, $y = (y_a)$, and $z = (z_r)$, which evolve jointly according to the replicator dynamics specified in Section~\ref{subsec:replicator}.

\subsection{Demand and Market Price}\label{subsec:demand}

Each commodity is sold in a single destination market whose clearing price depends on the quantity delivered. Following Assumption~2, the demand for the commodity $i$ is modeled by a constant-elasticity (isoelastic) function:
\begin{equation}
Q_i(S_i) = A_i\, S_i^{-\varepsilon_i}, 
\qquad A_i, \varepsilon_i > 0,
\label{eq:demand}
\end{equation}
where $S_i$ is the market price, $\varepsilon_i$ is the constant price elasticity of demand, and $A_i$ is a scale parameter that governs the overall size of the market. Inverting \eqref{eq:demand} gives the inverse-demand relationship:
\begin{equation}
S_i(Q_i) = \kappa_i\, Q_i^{-1/\varepsilon_i},
\qquad \kappa_i = A_i^{1/\varepsilon_i},
\label{eq:inverse-demand}
\end{equation}
which expresses the market-clearing price as a decreasing function of the quantity supplied to the market. The isoelastic form is adopted because the elasticity $\varepsilon_i$ admits a direct economic interpretation: a one percent increase in price reduces demand by $\varepsilon_i$ percent regardless of the price level, and because it guarantees a well-behaved, strictly downward-sloping demand curve for all positive prices and quantities.

\paragraph{Calibration.}
The scale parameter $A_i$ is not estimated independently, but is calibrated so that the demand curve passes through the observed nominal price--quantity pair. Given a benchmark price $S_i$ and the associated available volume $D_i$, imposing $Q_i(S_i) = D_i$ on \eqref{eq:demand} yields the following results.
\begin{equation}
A_i = D_i\, S_i^{\,\varepsilon_i},
\label{eq:calibration}
\end{equation}
so that the demand intercept is fully determined once the elasticity $\varepsilon_i$ is fixed. The elasticities themselves are taken from the empirical literature (Section~\ref{sec:case}), reflecting the structural price responsiveness of each commodity.

\paragraph{Role in the model.}
Consistent with Assumption~10, the market price $S_i$ is exogenous to any individual supplier's bidding decision, so the supplier and carrier payoffs of Section~\ref{subsec:payoffs} employ the benchmark price directly. The worst-case availability bound $\underline{D}_i = \varphi D_i$ used throughout the evolutionary stage is the one derived in 
Section~\ref{subsec:aro-counterpart}.

\subsection{Feasibility and Shared Constraints}\label{subsec:constraints}

Regardless of the strategies adopted by the three players, every shipment profile must satisfy a common set of physical and regulatory constraints. These constraints couple the players' decisions and define the feasible region within which the carrier's allocation is determined.

\paragraph{Capacity.}
The total volume shipped across all commodities cannot exceed the carrier's aggregate annual capacity:
\begin{equation}
\sum_{i \in \mathcal{I}} u_i \;\le\; T.
\label{eq:capacity}
\end{equation}

\paragraph{Availability.}
The volume shipped of each commodity cannot exceed the volume available for shipment. Under the worst-case realization of the uncertainty set, the relevant upper bound is the worst-case availability $\underline{D}_i = \varphi D_i$ introduced in Section~\ref{subsec:demand}:
\begin{equation}
0 \;\le\; u_i \;\le\; \underline{D}_i 
\qquad \forall\, i \in \mathcal{I}.
\label{eq:availability}
\end{equation}

\paragraph{Quota feasibility.}
The quota vector selected by the regulator must itself be feasible, in the sense that the mandated volumes are individually attainable and do not in aggregate exceed capacity:
\begin{equation}
0 \;\le\; L_i \;\le\; \underline{D}_i 
\quad \forall\, i \in \mathcal{I},
\qquad 
\sum_{i \in \mathcal{I}} L_i \;\le\; T.
\label{eq:quota-feasibility}
\end{equation}
Because the quota menu $\mathcal{R}$ is generated by the ARO stage subject to these same restrictions, every policy $L^{(r)} \in \mathcal{R}$ satisfies 
\eqref{eq:quota-feasibility} by construction.

\paragraph{Quota enforcement.}
Finally, the carrier must honor the mandated quota for every commodity, i.e., the shipped volume of each commodity is required to meet or exceed its quota.
\begin{equation}
u_i \;\ge\; L_i 
\qquad \forall\, i \in \mathcal{I}.
\label{eq:quota-enforcement}
\end{equation}
Constraint \eqref{eq:quota-enforcement} is the mechanism through which regulation binds: it forces the carrier to transport socially valuable commodities that it would otherwise exclude, and this hard lower bound is what distinguishes the regulated allocation from the unregulated, purely profit-driven one.

Constraints \eqref{eq:capacity}--\eqref{eq:quota-enforcement} together define the bounded feasible region over which the carrier optimizes. The manner in which the carrier allocates capacity within this region, first honoring the quota lower bounds \eqref{eq:quota-enforcement}, then distributing the residual capacity subject to \eqref{eq:capacity} and \eqref{eq:availability}, is described in Section~\ref{subsec:allocation}.

\subsection{Payoff Functions}\label{subsec:payoffs}

Each player evaluates its strategies according to a payoff function that maps a shipment profile to a scalar value. The three payoff functions are defined below.

\paragraph{Supplier profit.}
For a given shipped volume $u_i$ and a posted bid $p_i$, the supplier's profit $i$ comprises the revenue from the marketed volume, net of the transport payment to the carrier, less an inventory penalty on any unshipped volume:
\begin{equation}
\pi_i(u_i, p_i) 
= u_i \bigl( S_i - p_i \bigr) 
- g_i \bigl( D_i - u_i \bigr),
\label{eq:supplier-payoff}
\end{equation}
where $S_i$ is the exogenous market price, $g_i$ is the inventory/backlog penalty per unit, and $D_i$ is the volume the supplier has available. The first term is the net margin earned on the volume shipped; the second term is the penalty incurred on the residual $D_i - u_i$ that cannot be transported. A supplier thus faces a trade-off: bidding higher improves its priority in the carrier's allocation and reduces unshipped volume, but lowers the net margin $S_i - p_i$ earned on each shipped unit.

\paragraph{Carrier profit and the MRE cap.}
The carrier earns the transport margin $p_i - c_i$ on each unit shipped, where $c_i$ is the operating cost per unit. However, revenue in each flow is capped by the Maximum Revenue Entitlement $\mathrm{MRE}_i$, according to Assumption~5. The carrier's total profit is therefore:
\begin{equation}
\Pi^{\mathrm{car}}(u, p) 
= \sum_{i \in \mathcal{I}} 
\min \bigl\{ (p_i - c_i)\, u_i,\; \mathrm{MRE}_i \bigr\}.
\label{eq:carrier-payoff}
\end{equation}
The cap binds when the transport revenue in a flow would otherwise exceed $\mathrm{MRE}_i$, and prevents the carrier from extracting excessive revenue from any single commodity. Since $\mathrm{MRE}_i$ is calibrated to the maximum revenue attainable in full shipment (Section~\ref{sec:case}), the cap is consistent with the feasible bid range in \eqref{eq:price-grid}.

\paragraph{Regulator utility.}
The regulator values the social benefit of transporting each commodity, weighted by its social importance, net of the cost of enforcing the mandated quotas:
\begin{equation}
U^{\mathrm{reg}}(u, L) 
= \sum_{i \in \mathcal{I}} \beta_i\, u_i 
- \sum_{i \in \mathcal{I}} K_i(L_i),
\label{eq:regulator-payoff}
\end{equation}
where $\beta_i$ is the social welfare weight of commodity $i$ and $K_i(L_i)$ is the separable enforcement cost of Eq.~\eqref{eq:K-def}.

The first term of Eq. \eqref{eq:regulator-payoff} reflects the value to society of the actual delivered goods, whereas the second term captures the total cost of enforcement. Using the function $K_i(L_i)$, the regulator internalizes the enforcement cost. In other words, raising the quota increases the value to society from the guaranteed deliveries, but also raises both the linear and convex part of the enforcement cost. The interaction leads to the Pareto frontier identified at the ARO stage.

\subsection{The Carrier's Allocation Procedure}\label{subsec:allocation}

For a quota vector $L$, a set of bids by suppliers $p$, and an allocation scheme $a \in \mathcal{A}$, the carrier computes the shipment vector $u$ using a two-stage process that implements the constraints discussed in Section~\ref{subsec:constraints}. This process is realized through the subroutine $\textsc{Allocate}(L, p, a)$, called by Algorithm~\ref{alg:aro-egt}.

\paragraph{Step 1: Honour mandated quotas.}
The carrier first assigns each commodity its mandated quota, subject to worst-case availability:
\begin{equation}
u_i \;\gets\; \min\bigl\{ L_i,\; \underline{D}_i \bigr\}
\qquad \forall\, i \in \mathcal{I}.
\label{eq:alloc-step1}
\end{equation}
This step guarantees that the quota-enforcement constraint \eqref{eq:quota-enforcement} is satisfied before any profit-driven allocation takes place.

\paragraph{Step 2: Distribute residual capacity.}
The residual capacity remaining after the quotas have been honored, $T_{\mathrm{res}} = T - \sum_{i} u_i$, is then distributed among the commodities according to the chosen policy $a$, up to the remaining availability of each commodity $\underline{D}_i - u_i$. The order and manner in which the residual capacity is assigned depend on the policy, as specified below.

The set $\mathcal{A}$ comprises five behaviorally distinct allocation policies, chosen to span the range from purely profit-oriented behavior to socially oriented behavior:

\begin{itemize}
    \item \textbf{$a_0$ (Pure profit).} Residual capacity is assigned to commodities in decreasing order of transport margin $p_i - c_i$. Commodities with negative margin do not receive residual capacity beyond their mandated quota, reflecting that a profit-maximizing carrier does not voluntarily ship at a loss.

    \item \textbf{$a_1$ (Social first).} Residual capacity is assigned in decreasing order of social weight $\beta_i$, prioritizing socially valuable commodities irrespective of their transport margin.

    \item \textbf{$a_2$ (Proportional).} Residual capacity is distributed among commodities in proportion to their remaining availability $\underline{D}_i - u_i$.

    \item \textbf{$a_3$ (MRE-greedy).} Residual capacity is assigned in decreasing order of effective margin, where each commodity's per-unit  margin is capped at $\mathrm{MRE}_i / \underline{D}_i$. Commodities whose effective margin is negative receive no residual capacity beyond their quota.

    \item \textbf{$a_4$ (Balanced).} Residual capacity is assigned in decreasing order of a composite score that weights normalized transport margin and normalized social value equally, representing a carrier that balances private profit against social considerations.
\end{itemize}

In every policy, residual assignment proceeds until either the residual capacity is exhausted or no commodity retains unfilled availability. The distinction between profit-oriented policies ($a_0$, $a_3$), which decline to ship negative-margin commodities beyond their quota, and socially oriented policies ($a_1$, $a_4$), which do not, is central to the evolutionary dynamics. This distinction determines whether a socially valuable but loss-making commodity reaches the market through voluntary allocation or only through the mandated quota. This, in turn, governs the relative fitness of the quota policies of the regulator, since a commodity that the carrier will not ship voluntarily can be delivered only by binding the enforcement constraint \eqref{eq:quota-enforcement}.

\subsection{Replicator Dynamics}\label{subsec:replicator}

The three populations adapt their strategies over discrete time according to replicator dynamics with mutation. At each step, the expected payoff of each strategy is evaluated against the prevailing behavior of the other two populations, and the strategy shares are updated so that better-performing strategies gain ground. The two components of this process, payoff evaluation and share updating, are described below.

\paragraph{Mean-field payoff evaluation.}
Because each player faces a population of opponents rather than a single opponent, payoffs are evaluated in the mean-field sense. Hence, each strategy is assessed against the \emph{average} behavior of the other two populations, weighted by their current strategy shares. Let $\bar{p}_i = \sum_{j \in J_i} x_{ij}\, p_{ij}$ denote the mean bid of supplier $i$. The expected payoff of carrier policy $a$ and of regulator policy $r$ is obtained by evaluating the allocation $\textsc{Allocate}(L^{(r)}, \bar{p}, a)$ between opposing populations:
\begin{align}
\Pi^{\mathrm{car}}_a 
&= \sum_{r=1}^{R} z_r\, 
\Pi^{\mathrm{car}}\!\bigl( \textsc{Allocate}(L^{(r)}, \bar{p}, a),\, 
\bar{p} \bigr),
\label{eq:mf-carrier}\\[2pt]
U^{\mathrm{reg}}_r 
&= \sum_{a \in \mathcal{A}} y_a\, 
U^{\mathrm{reg}}\!\bigl( \textsc{Allocate}(L^{(r)}, \bar{p}, a),\, 
L^{(r)} \bigr).
\label{eq:mf-regulator}
\end{align}
The expected payoff of the supplier-$i$ strategy $p_{ij}$ is evaluated by allowing that supplier to deviate to $p_{ij}$ while the remaining suppliers bid at their mean and averaging the resulting profit over the carrier and regulator populations:
\begin{equation}
\pi_{ij} 
= \sum_{r=1}^{R} \sum_{a \in \mathcal{A}} 
z_r\, y_a\; 
\pi_i\!\bigl( u_i^{(r,a,j)},\, p_{ij} \bigr),
\label{eq:mf-supplier}
\end{equation}
where $u_i^{(r,a,j)}$ is the volume that commodity $i$ receives under $\textsc{Allocate}$ when supplier $i$ bids $p_{ij}$ and all other suppliers bid at their means.

\paragraph{Replicator update with mutation.}
Given the expected payoffs, each population share is revised using the discrete-time replicator rule augmented with a small mutation term. For a population with strategy shares $w = (w_k)$ and payoffs $(\phi_k)$, the update is:
\begin{equation}
w_k' = (1 - \mu)\, \frac{w_k\, \tilde{\phi}_k}{\sum_{\ell} w_\ell\, \tilde{\phi}_\ell} + \frac{\mu}{|\mathcal{K}|},
\label{eq:replicator}
\end{equation}
where $\tilde{\phi}_k = \phi_k - \min_{\ell} \phi_\ell + \delta$ is the payoff shifted to be strictly positive (with $\delta > 0$ a small constant), $|\mathcal{K}|$ is the number of strategies in the population and $\mu \in [0,1]$ is the mutation rate. This first term represents the traditional replicator dynamics in which strategies producing payoffs that are greater than the average will grow more common, whereas strategies producing below average payoffs will become less common. The second term allows for a small but constant mutation that keeps diversity within strategies alive. Applying \eqref{eq:replicator} to the supplier, carrier, and regulator populations, using the payoffs \eqref{eq:mf-carrier}--\eqref{eq:mf-supplier}, yields the coupled update:
\begin{equation}
x_{ij} \to x_{ij}', 
\qquad 
y_a \to y_a', 
\qquad 
z_r \to z_r',
\label{eq:coupled-update}
\end{equation}
which is iterated from uniform initial shares until the strategy distributions stabilize or a prescribed horizon $T_{\max}$ is reached. Because the three populations enter one another's payoffs through \eqref{eq:mf-carrier}--\eqref{eq:mf-supplier}, the update \eqref{eq:coupled-update} constitutes a genuine co-evolutionary system in which the fitness of each player's strategies depends on the concurrent state of the other two.

\subsection{Evolutionary Equilibrium}\label{subsec:equilibrium}

The iteration of the coupled update \eqref{eq:coupled-update} generates the trajectories of the strategy distributions $x(t)$, $y(t)$, and $z(t)$. The long-term behavior of these trajectories characterizes the outcome of the game.

\paragraph{Definition.}
A state $(x^{\ast}, y^{\ast}, z^{\ast})$ is an evolutionary equilibrium of the system if it is a rest point of the replicator dynamics, that is if applying \eqref{eq:replicator} to each population leaves the strategies' shares unchanged up to the mutation floor. In such a state, every strategy retained with positive share earns the population-average payoff, so that no strategy can increase its share through the replicator mechanism. In discrete implementation, the system is deemed to have reached equilibrium when the strategy distributions cease to change appreciably between successive iterations or when the prescribed horizon $T_{\max}$ is attained.

\paragraph{Possible outcomes.}
The replicator dynamics admits three qualitatively distinct long-run outcomes, each with a different interpretation for the regulated market:
\begin{enumerate}
    \item \emph{Pure (dominant) equilibrium.} A single strategy in each population absorbs essentially all of the share, so that the market settles on a well-defined and predictable configuration of quota policy, bids, and allocation rule.
    \item \emph{Stable mixed equilibrium.} Two or more strategies coexist at fixed positive shares, sustained by the fact that their relative payoffs  depend on the prevailing strategy mix. The market settles on a stable but diversified pattern of behavior rather than a single strategy.
    \item \emph{Persistent cycling.} The strategy shares fail to converge and instead oscillate indefinitely, reflecting a configuration in which no strategy profile is stable against the adaptation of the others.
\end{enumerate}
The outcome obtained is determined by the payoff structure and hence by the underlying commodity data, the enforcement-cost parameters, and the capacity and social-weight settings. The classification is therefore an empirical property of a given instance, not a feature imposed a priori, and it is examined for the Canadian case in Section~\ref{sec:results}.

\paragraph{Recovery of the equilibrium outcome.}
Once the dynamics have stabilized, the prevailing strategy of each population is identified by its dominant share, where the equilibrium carrier policy is $a^{\ast} = \arg\max_{a} y_a^{\ast}$, the equilibrium regulator policy is $r^{\ast} = \arg\max_{r} z_r^{\ast}$ with the associated quota vector $L^{\ast} = L^{(r^{\ast})}$, and the equilibrium bid of each supplier is $p_i^{\ast} = p_{i j^{\ast}}$ with $j^{\ast} = \arg\max_{j} x_{ij}^{\ast}$. The equilibrium shipment profile is then recovered by evaluating the carrier's allocation procedure at these dominant strategies,
\begin{equation}
u^{\ast} = \textsc{Allocate}\bigl( L^{\ast}, p^{\ast}, a^{\ast} \bigr),
\label{eq:equilibrium-allocation}
\end{equation}
from which the equilibrium payoffs of all three players, and the realized social welfare $\sum_{i} \beta_i u_i^{\ast}$, follow directly. This equilibrium outcome $(L^{\ast}, p^{\ast}, a^{\ast}, u^{\ast})$ constitutes the solution of the integrated framework and forms the basis for the results reported in Section~\ref{sec:results}.

\section{Case Study: Canadian Rail Freight}\label{sec:case}

The proposed framework is applied to a case of rail freight transport in Canada. First, data sources and units of measure are described and then the baseline parameters and calibration of the derived parameters are provided.

\subsection{Data Sources and Unit Conventions}
\label{subsec:data}

The four commodities considered as part of the case study are representative of the rail freight commodities that caused the capacity crisis in 2013–2014. The four commodities include crude oil, corn, barley, and oat. In essence, those commodities represent the conflict that lies behind the need for the ``Fair Rail for Grain Farmers Act,'' since there is the question of a highly profitable energy commodity (crude oil) and less profitable grain commodities. Parameter values are obtained from the sources used by \citet{pazoki2024}, making the results comparable to the deterministic benchmark.

All values are presented on an annual basis, measured in identical units for each commodity. Commodity volumes are presented in bushels per year using the weight conversion for each commodity, which is 56~lb for corn, 48~lb for barley, and 32~lb for oat. However, a bushel is a unit of volume, not of mass. Therefore, the volume of crude oil will be measured in bushels on the basis of the conversion factor of one cubic meter being equal to $28.38$ bushels. Values are measured in Canadian dollars (CAD). Production volumes refer to the 2019 harvest season, while per unit costs are based on provincial cost surveys for 2020. Crude oil prices are defined on the basis of the 2019 price of the Mixed Sweet Blend (Edmonton) and using the 2020 volume forecast of crude oil production.

Because aggregate rail capacity cannot be fully realized when only four commodities are considered, the capacity parameter $T$ is treated as a design variable and examined in a range of values in the sensitivity analysis of Section~\ref{sec:results}, rather than being fixed at a single system-wide figure. The baseline value adopted for the primary analysis is $T = 650{,}000{,}000$ bushels per year, according to the benchmark study.

\subsection{Baseline Parameters}\label{subsec:baseline-parameters}

Table~\ref{tab:baseline} reports the baseline product-level parameters used in the payoff functions, namely, the annual available volume $D_i$, the exogenous market price $S_i$, the carrier's per-unit operating cost $c_i$, and the supplier's per-unit inventory penalty $g_i$. The last column reports the maximum transport margin per unit by the carrier $S_i + g_i - c_i$, which determines whether a commodity is shipped voluntarily by a profit-oriented carrier.

\begin{table}[H]
\centering
\caption{Baseline product parameter values.}
\label{tab:baseline}
\renewcommand{\arraystretch}{1.25}
\begin{tabular}{lccccc}
\toprule
\textbf{Commodity} & \textbf{$D_i$ (bu/yr)} & \textbf{$S_i$ (CAD)} 
& \textbf{$c_i$ (CAD)} & \textbf{$g_i$ (CAD)} 
& \textbf{$m_i$ (Eq.~\ref{eq:margin})} \\
\midrule
Crude oil & $128{,}185{,}505$ & $15.3000$ & $1.2970$ & $0.7900$ & $+14.7930$ \\
Corn      & $440{,}916{,}666$ & $0.0900$  & $0.0600$ & $0.0042$ & $+0.0342$ \\
Barley    & $440{,}924{,}524$ & $0.1110$  & $0.1130$ & $0.0049$ & $+0.0029$ \\
Oat       & $275{,}577{,}827$ & $0.1100$  & $0.1180$ & $0.0073$ & $-0.0007$ \\
\bottomrule
\end{tabular}
\end{table}

The specific values of the parameters explain the basic asymmetry that lies at the heart of the problem under regulation. Crude oil has an enormous transport margin of $14.79$ CAD per bushel, several orders of magnitude higher than any of the grains, which makes the profit-maximizing shipper move crude oil first and only then move anything else. Corn and barley have small but non-negative margins ($0.0342$ and $0.0029$ CAD per bushel, respectively). Therefore, they are shipped only after all crude oil has been shipped, and to the extent that there is capacity, only because the carrier has no choice. However, oat has a negative maximum margin ($-0.0007$ CAD per bushel). Thus, no matter how high the bid, the shipper will not even receive a marginal transport margin by shipping this commodity. Hence, oat, which is socially valuable but unprofitable to transport, is precisely the commodity that should be protected by regulation. It is the absolute smallness of the margin and, even more so, its sign than its value that matters. Any commodity with a negative maximum transport margin will be out of the picture in the voluntary allocation scheme. The transport cost sensitivity studied in Section~\ref{subsec:sensitivity} addresses this issue directly and shows that the structure of equilibrium depends on the sign of the margin. Hence, this case study provides an illustration of the general problem formalized in Sections~\ref{sec:problem} and \ref{sec:egt}.

\subsection{Calibration of Derived Quantities}
\label{subsec:calibration}

In addition to the primitive parameters listed in Table~\ref{tab:baseline}, some additional variables must be specified to complete the specification of the model, namely, demand elasticities and intercepts, maximum revenue entitlements, enforcement cost parameters, social welfare weights, and the worst-case availability factor. The calibration of all these variables is provided below and summarized in Table~\ref{tab:derived}.

\paragraph{Demand elasticities and intercepts.}
The price elasticities are taken from empirical studies on commodity demand. A short-term elasticity of $\varepsilon = 0.30$ is chosen for crude oil, capturing the well-known price inelasticity of petroleum demand \citep{hamilton2009}. The elasticities of $0.50$ for corn and $0.60$ for barley and oat are assumed based on the elasticities of staple agricultural goods \citep{roberts2013}. Based on these elasticities, the demand intercept for each good is set using Eq. \eqref{eq:calibration} to ensure that the isoelastic demand curve generates the observed nominal quantity at the benchmark price: $A_i = D_i\, S_i^{\varepsilon_i}$.

\paragraph{Maximum revenue entitlement.}
The MRE ceiling of each flow is defined as the maximum of the transport revenue generated under total shipment of the commodity at the highest rational bid price, i.e., $\mathrm{MRE}_i = (S_i + g_i - c_i)\, D_i$. This calibration guaranties that the constraint is compatible with the feasible bid range $[0, S_i + g_i]$ of Eq. \eqref{eq:price-grid}, where the supplier placing the highest possible bid and shipping the entire available volume will exactly generate the entitlement $\mathrm{MRE}_i$.

\paragraph{Enforcement-cost parameters.}
The enforcement-cost function \eqref{eq:K-def} is calibrated so that the enforcement cost appropriately depends on the size of the commodity market. The convex curvature is taken inversely proportional to the available volume, $\eta_i = \eta_0 / D_i$ with $\eta_0 = 0.0293$, making the enforcement cost of imposing the same fraction of the available volume the same for all commodities regardless of their size. The linear component is calibrated using the same CAD per unit rate $\lambda_i = 0.006$ per bushel, and the fixed enforcement cost $\rho_i$ is chosen to be negligible compared to the linear and convex parts in the baseline analysis. As a result, the maximum useful enforcement expenditure is $K^{\max} = 4.76$ million CAD, which is the amount needed to enforce the total shipment of the worst-case availability of oat, as shown in Section~\ref{subsec:aro-pareto}.

\paragraph{Social welfare weights.}
In this study, we use the unequally important case of \citet{pazoki2024}, where all commodities except oat have a weight of $\beta = 0.1$, while oat have a much higher weight of $\beta_{\mathrm{oat}} = 1.0$. This corresponds to the realistic situation where a non-commercially valuable commodity is socially important, creating a tension between the private profit and social welfare because of the food security concerns. The sensitivity of the results to $\beta_{\mathrm{oat}}$ is tested in Section \ref{sec:results}.

\paragraph{Worst-case availability factor.}
The uncertainty of commodity availability is captured by a worst-case availability factor $\varphi = 0.9$ in the baseline analysis, leading to a downward deviation from the nominal availability of $\underline{D}_i = 0.9\, D_i$. The effect of more pronounced uncertainty and even severe supply disruptions is tested by changing $\varphi$ in the sensitivity analysis section.

\begin{table}[H]
\centering
\footnotesize
\caption{Calibrated derived quantities.}
\label{tab:derived}
\renewcommand{\arraystretch}{1.2}
\setlength{\tabcolsep}{4pt}
\begin{tabularx}{\textwidth}{l c c X}
\toprule
\textbf{Quantity} & \textbf{Crude / Corn / Barley / Oat} & 
\textbf{Value / Rule} & \textbf{Source} \\
\midrule
Elasticity $\varepsilon_i$ & $0.30 / 0.50 / 0.60 / 0.60$ & literature & 
\citep{hamilton2009, roberts2013} \\
Demand intercept $A_i$ & $D_i\, S_i^{\varepsilon_i}$ & calibrated; sets $\varphi$ & 
Eq.~\eqref{eq:calibration} \\
$\mathrm{MRE}_i$ & $(S_i + g_i - c_i)\, D_i$ & calibrated & Assumption~5 \\
Convex cost $\eta_i$ & $\eta_0 / D_i,\ \eta_0 = 0.0293$ & scaled & --- \\
Linear cost $\lambda_i$ & $0.006$ CAD/bu & common & --- \\
Social weight $\beta_i$ & $0.1 / 0.1 / 0.1 / 1.0$ & scenario & 
\citep{pazoki2024} \\
Worst-case factor $\varphi$ & $0.9$ (baseline) & assumed & --- \\
\bottomrule
\end{tabularx}
\end{table}

\section{Results and Discussion}
\label{sec:results}

The results of the application of the proposed framework to the Canadian rail example follow in this section. First, the Pareto front obtained using the ARO phase is provided, followed by evolutionary equilibrium, the market without regulation, and the sensitivity analysis. All results are obtained with a frontier of $R = 20$ Pareto policies, supplier price grids of $|J_i| = 10$ points, a mutation rate of $\mu = 0.005$, and a replicator horizon of $T_{\max} = 50$, using the calibrated parameters of Section~\ref{sec:case}. The strategy distributions stabilize in roughly forty iterations. Longer horizons were also examined and yielded the same equilibrium. Consequently, the shorter horizon is reported as it displays the transient dynamics more legibly.

\subsection{The Robust Pareto Frontier}
\label{subsec:pareto-results}

\begin{figure}[H]
\centering
\includegraphics[width=\textwidth]{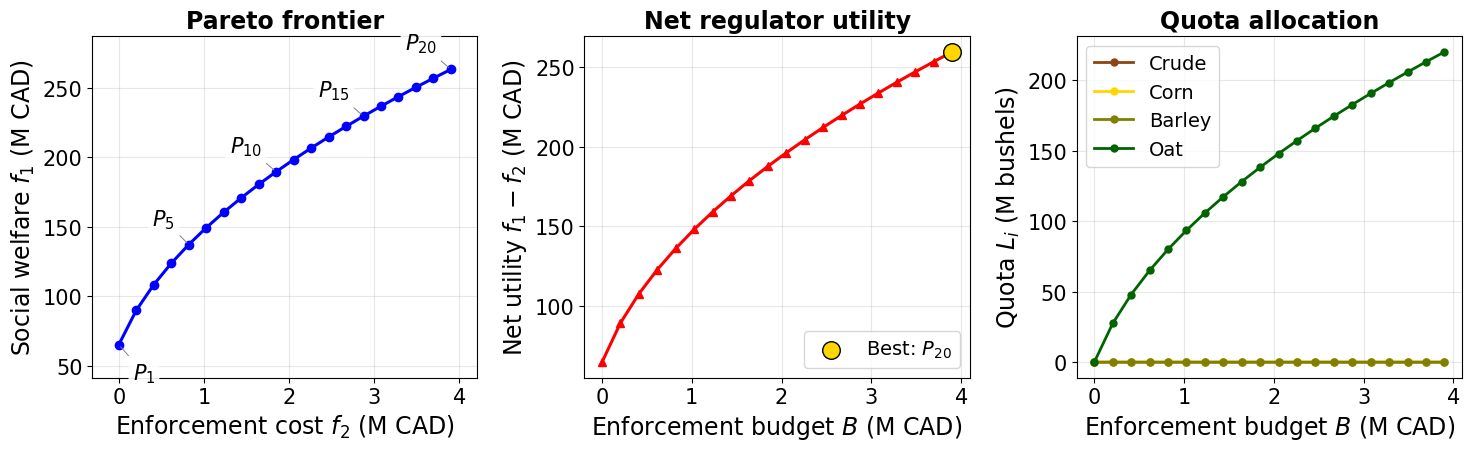}
\caption{The robust Pareto frontier. \\Left: social welfare $f_1$ against enforcement cost $f_2$ across the twenty Pareto policies. \\Center: net regulator utility $f_1 - f_2$ as a function of the enforcement budget. \\Right: the mandated quota of each commodity as the budget increases.}
\label{fig:pareto}
\end{figure}

The ARO stage generates a frontier of 20 robust and Pareto-efficient quota policies by solving the $\varepsilon$-constraint problem across enforcement budgets ranging from zero to $K^{\max} = 4.76$ million CAD. This upper limit is the expenditure required to protect the full worst-case availability of every commodity that the carrier will not ship voluntarily. Since oat is the only such commodity, $K^{\max} = K_{\mathrm{oat}}(\underline{D}_{\mathrm{oat}})$ with $\underline{D}_{\mathrm{oat}} = 248.0$ million bushels. Beyond this point, an additional budget does not bring more benefits because every remaining commodity moves without a mandate. Figure~\ref{fig:pareto} displays the resulting trade-off between worst-case social welfare ($f_1 = \sum_i \beta_i u_i$) and enforcement cost ($f_2 = \sum_i K_i(L_i)$), together with net regulator utility $f_1 - f_2$ and associated quota allocations.

Table~\ref{tab:pareto-summary} reports five representative policies that span the frontier, from no intervention ($P_1$) to full protection ($P_{20}$).

\begin{table}[H]
\centering
\footnotesize
\caption{Representative policies along the robust Pareto frontier. All quotas other than oat are zero across the entire frontier and are therefore omitted.}
\label{tab:pareto-summary}
\renewcommand{\arraystretch}{1.3}
\setlength{\tabcolsep}{5pt}
\begin{tabularx}{\textwidth}{l Y Y Y Y Y}
\toprule
\textbf{Policy} & \textbf{$B_r$ (M CAD)} & \textbf{$f_1$ welfare (M CAD)} 
& \textbf{$f_2$ enforce. (M CAD)} & \textbf{$f_1\!-\!f_2$ (M CAD)} 
& \textbf{Oat quota (M bu)} \\
\midrule
$P_{1}$   & $0.00$ & $65.0$  & $0.00$ & $65.0$  & $0.0$   \\
$P_{5}$   & $1.00$ & $147.8$ & $1.00$ & $146.8$ & $92.0$  \\
$P_{10}$  & $2.25$ & $206.4$ & $2.25$ & $204.1$ & $157.1$ \\
$P_{15}$  & $3.51$ & $250.9$ & $3.51$ & $247.4$ & $206.5$ \\
$P_{20}$  & $4.76$ & $288.2$ & $4.76$ & $283.5$ & $248.0$ \\
\bottomrule
\end{tabularx}
\end{table}

There are three distinctive properties of the frontier. First, social welfare rises monotonically along with the enforcement budget, starting with $65.0$ million CAD when the budget is zero and reaching $288.2$ million CAD at the full protection level and showing an increase of more than four times compared to the frontier’s null policy. Second, the frontier is concave, reflecting the convex nature of the enforcement cost function. In particular, the marginal benefit of enforcement decreases from $150$ CAD of welfare for each dollar spent on enforcement at $L_{\mathrm{oat}} = 0$ to $27.8$ under full protection level. The marginal benefit has a closed form. Since a bushel of oat quota replaces a bushel of corn, the marginal benefit is equal to $\beta_{\mathrm{oat}} - \beta_{\mathrm{corn}} = 0.9$ and the marginal cost is equal to $\lambda_{\mathrm{oat}} + \eta_{\mathrm{oat}} L_{\mathrm{oat}}$. Therefore, 
\begin{equation}
\frac{d f_1}{d B} 
= \frac{\beta_{\mathrm{oat}} - \beta_{\mathrm{corn}}}
       {\lambda_{\mathrm{oat}} + \eta_{\mathrm{oat}} L_{\mathrm{oat}}}.
\label{eq:shadow-budget}
\end{equation}
The frontier is bounded not because of prohibitively high enforcement cost, but because oat is exhausted, and the last dollar of enforcement still generates welfare more than twenty-eight times higher than the cost.

Third, the entire enforcement budget is spent on one commodity. As the right panel of Figure~\ref{fig:pareto} shows, the mandated quotas for crude oil, corn, and barley remain zero across the whole frontier, while the oat quota increases steadily with the budget. This follows from the margin structure of Table~\ref{tab:baseline}, i.e., crude oil, corn, and barley all have positive $m_i$ and are shipped voluntarily by a profit-oriented carrier, while oat, the only commodity with a negative $m_i$, reaches the market only if the regulator mandates it.

The frontier constitutes the regulator's strategy menu $\mathcal{R}$ for the evolutionary stage, ranging from the null policy $P_1$ to the full-protection policy $P_{20}$.

\subsection{Evolutionary Equilibrium}
\label{subsec:egt-results}

\begin{figure}[H]
\centering
\includegraphics[width=\textwidth]{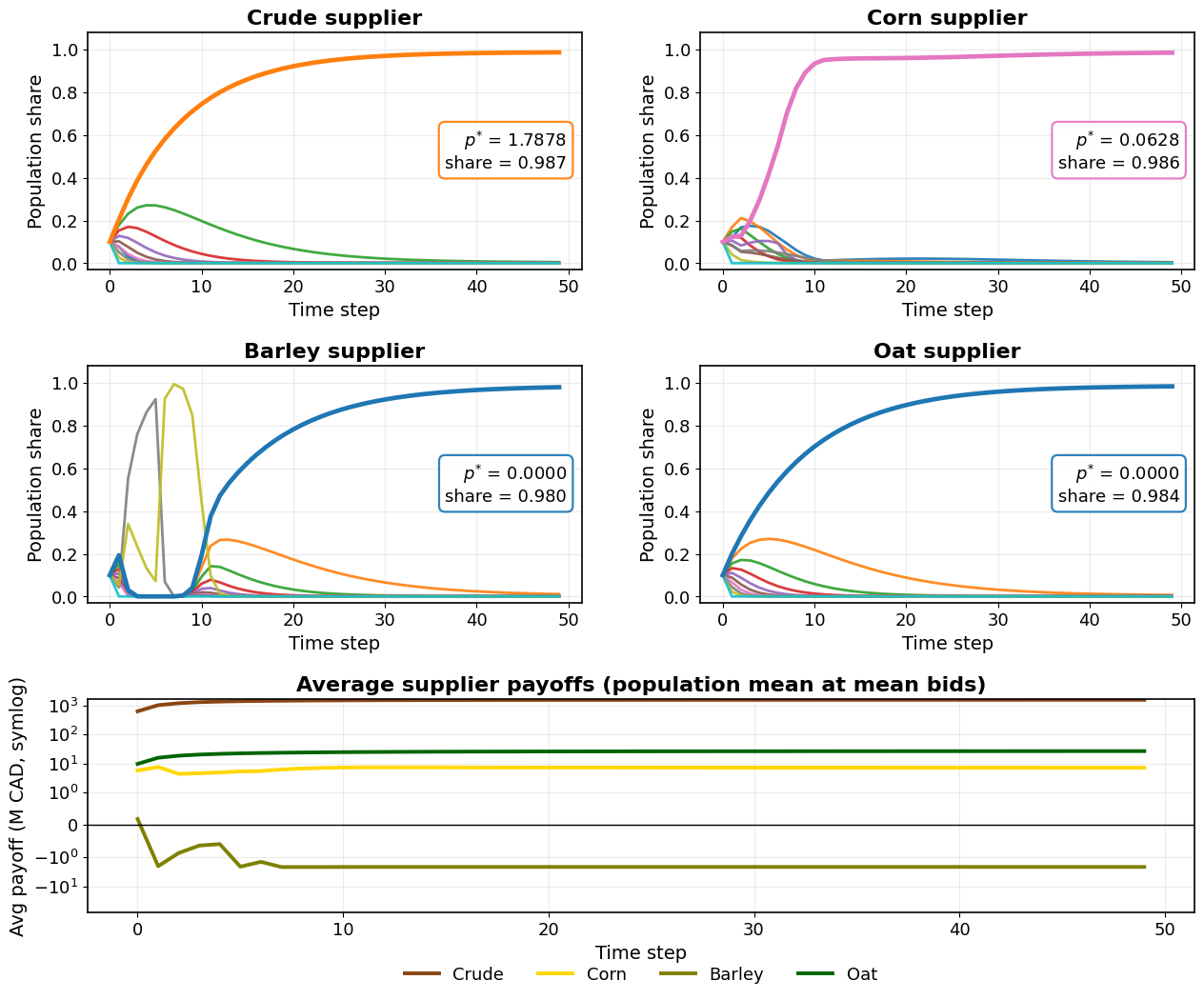}
\caption{Evolution of the four supplier populations.\\Each of the upper panels shows the distribution over the ten price-bid strategies of one commodity, with the dominant strategy drawn in bold. \\The lower panel shows average supplier payoffs on a symmetric logarithmic scale, taken as population means evaluated at mean bids.}
\label{fig:egt-suppliers}
\end{figure}

\begin{figure}[H]
\centering
\includegraphics[width=\textwidth]{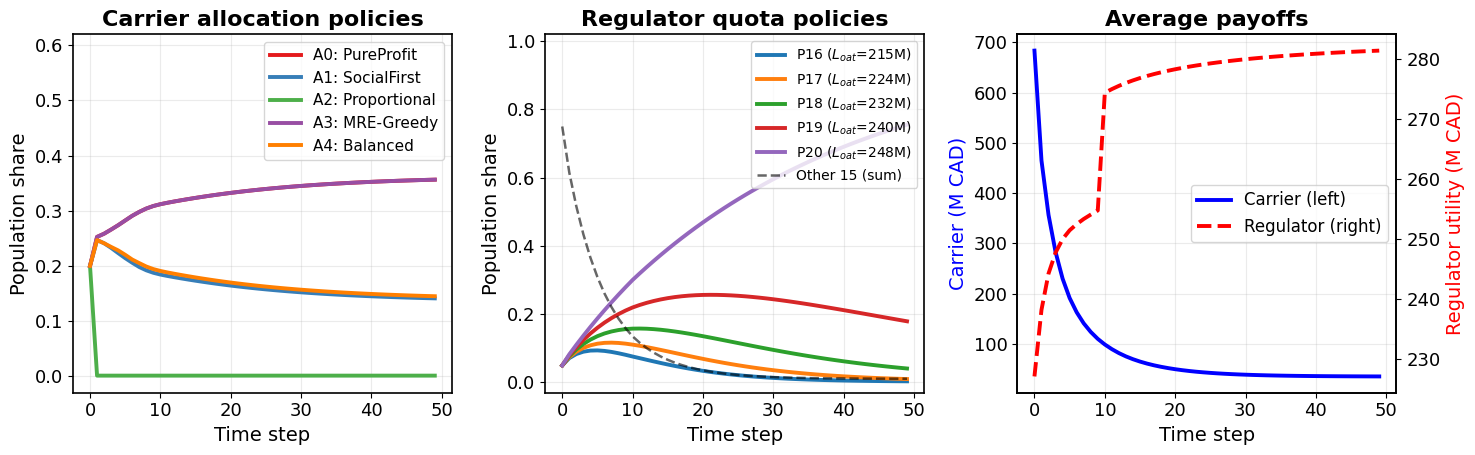}
\caption{Evolution of the carrier and regulator populations. \\Left: shares of the five allocation policies. \\Center: shares of the five largest quota policies, with the remaining fifteen aggregated as a dashed line. \\Right: average payoffs, taken as population means at mean bids rather than the equilibrium values at the dominant strategies reported in Table~\ref{tab:egt-equilibrium}.}
\label{fig:egt-carrier}
\end{figure}

With the Pareto frontier fixed as the regulator's strategy menu, the three populations co-evolve under the replicator dynamics of Section~\ref{subsec:replicator} from uniform initial shares. Figures~\ref{fig:egt-suppliers} and~\ref{fig:egt-carrier} trace the evolution of the supplier, carrier, and regulator strategy distributions, together with the associated payoff paths, over the simulation horizon.

In this model, convergence towards a nearly pure strategy is observed, as described in Table~\ref{tab:egt-equilibrium}. The regulator group focuses exclusively on the full-protection strategy $P_{20}$, which occupies a predominant position $0.762$. At the same time, neighboring strategies $P_{19}$ and $P_{18}$ occupy positions of $0.176$ and $0.039$, respectively. This arrangement is explained by the fact that the utility function $f_1 - f_2$ on the frontier increases monotonically, so the best-positioned policy from this perspective is the most evolutionary advantageous. Suppliers demonstrate a similar pronounced focus of their bids, with dominants from $98.1\%$ to $98.7\%$.

The carrier population divides into four policies rather than settling on one. The pure-profit policy $a_0$ and the MRE-greedy policy $a_3$ each achieve $0.357$, while the socially oriented policies $a_1$ and $a_4$ retain $0.141$ and $0.145$. Two mechanisms produce this pattern. The tie between $a_0$ and $a_3$ arises because the two differ only in that $a_3$ caps each commodity's per-unit margin at $\mathrm{MRE}_i / \underline{D}_i$ before ranking, so they diverge only when the revenue entitlement binds. It does not bind here since at equilibrium the carrier earns $56.6$ million CAD on crude oil against an entitlement of $1{,}896$ million, and $1.1$ million on corn against an entitlement of $15.1$ million, slack by factors of roughly 33 and 14. Equilibrium bids settle far below their ceilings, since a supplier bids only enough to secure priority, and shipments are bounded by $\underline{D}_i = 0.9 D_i$ rather than $D_i$. As in \citet{pazoki2024}, the entitlement is therefore not restrictive in this case study.

The survival of $a_1$ and $a_4$ in non-trivial shares reflects a second effect. Once oat is fully protected, only $286.6$ million bushels of residual capacity remain to be allocated, and all three commodities competing for capacity carry the same social weight. Disparities in payoffs between profit-maximizing and socially responsible allocation mechanisms become smaller, and there is less pressure on firms to choose the former. Therefore, quota management does more than simply redistribute. That is, by cutting off the controversial share from the carrier’s discretion, quotas reduce the impact that the choice of mechanism can have on the result.

\begin{table}[H]
\centering
\footnotesize
\caption{Evolutionary equilibrium: dominant strategies and equilibrium allocation.}
\label{tab:egt-equilibrium}
\renewcommand{\arraystretch}{1.3}
\setlength{\tabcolsep}{5pt}
\begin{tabularx}{\textwidth}{l Y Y Y Y Y}
\toprule
\textbf{Commodity} & \textbf{Bid $p_i^{\ast}$ (CAD)} & \textbf{Margin $p_i^{\ast}\!-\!c_i$} & \textbf{Quota $L_i^{\ast}$ (M bu)} & \textbf{Shipped $u_i^{\ast}$ (M bu)} & \textbf{Profit (M CAD)} \\
\midrule
Crude oil & $1.7878$ & $+0.4908$ & $0.0$   & $115.4$ & $1{,}548.7$ \\
Corn      & $0.0628$ & $+0.0028$ & $0.0$   & $286.6$ & $7.1$       \\
Barley    & $0.0000$ & $-0.1130$ & $0.0$   & $0.0$   & $-2.2$      \\
Oat       & $0.0000$ & $-0.1180$ & $248.0$ & $248.0$ & $27.1$      \\
\midrule
\multicolumn{6}{l}{\textit{Dominant carrier policies: $a_0$/$a_3$ (0.357 each);\quad 
dominant regulator policy: $P_{20}$ (0.762)}} \\
\bottomrule
\end{tabularx}
\end{table}

\begin{figure}[H]
\centering
\includegraphics[width=\textwidth]{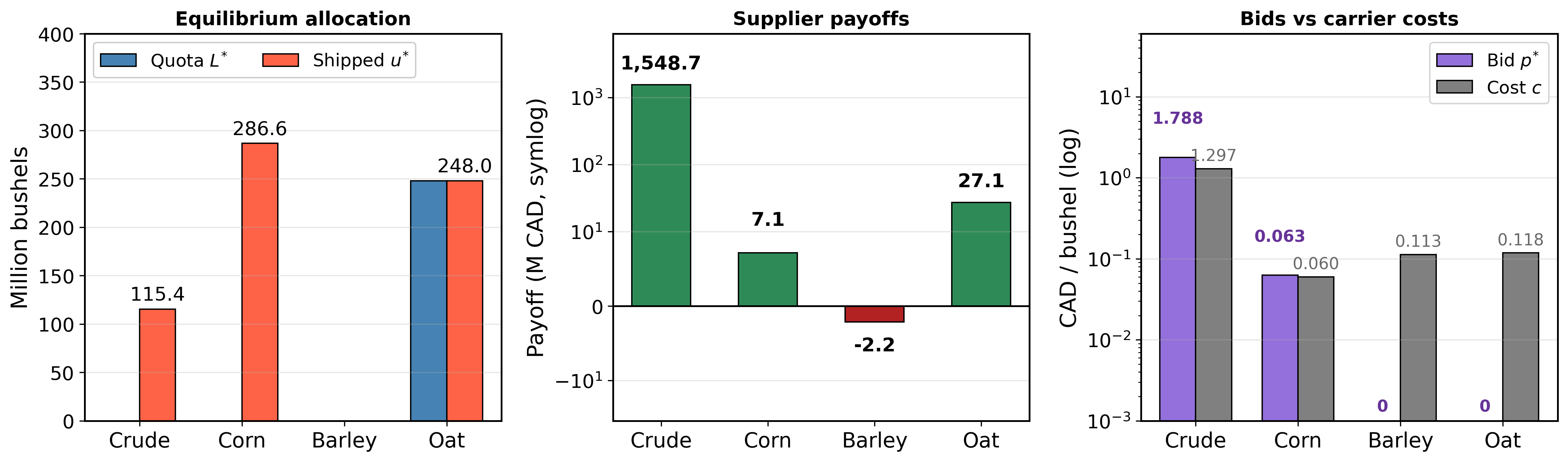}
\caption{Equilibrium outcome. \\Left: mandated quota against realized shipment for each commodity; oat alone ships exactly its mandate. \\Center: supplier payoffs on a symmetric logarithmic scale, with barley alone negative. \\Right: dominant bids against the carrier's per-unit operating cost on a logarithmic scale.}
\label{fig:equilibrium}
\end{figure}

The equilibrium bids reveal the strategic logic of the suppliers. Crude oil, whose transport margin is large and positive, bids $1.79$ CAD per bushel, the minimal bid sufficient to secure priority in the carrier's residual allocation, and ships its full available volume of $115.4$ million bushels. Corn bids just above the carrier's operating cost, retaining a margin of $0.0028$ CAD per bushel that keeps it attractive and absorbs the remaining capacity. Barley and oat both bid zero, but for different reasons. Oat cannot offer the carrier a profitable rate in any bid that a rational supplier would offer, since $m_{\mathrm{oat}} < 0$. Barley could be profitable in principle, since $m_{\mathrm{barley}} = +0.0029$, but the margin is so thin that the required bid would leave barley with almost no profit, and in any case it is dominated by crude oil and corn in the carrier's ranking. Therefore, both oat and barley maximize their surplus $S_i - p_i$ by not bidding and avoiding shipment.

The consequences of this choice differ sharply between the two. Oat, protected by a mandated quota of $248.0$ million bushels, is delivered in full and earns $27.1$ million CAD. Barley, which does not receive a quota, does not ship and records a loss of $2.2$ million CAD as it pays the inventory penalty of $g_i D_i$ on its entire available volume. Exclusion from the network is not merely a foregone opportunity for an unprotected supplier; it is a realized cost.

At equilibrium the carrier fills its entire capacity of $650$ million bushels, with crude oil, corn, and oat. The social welfare realized is $288.2$ million CAD, matching the welfare of policy $P_{20}$ on the Pareto frontier and confirming that the evolutionary process selects the upper extreme of the frontier as the stable outcome.

\subsection{Comparison with the Unregulated Market}
\label{subsec:unregulated}

The value of regulation is measured against a market with no quotas. In this benchmark the regulator takes no action and the carrier allocates capacity by transport margin. Suppliers keep their equilibrium bids, with one exception. At those bids barley would bid zero and $137.8$ million bushels of capacity would stay idle, yet any barley bid between the carrier's cost ($0.113$ CAD per bushel) and its maximum ($0.1159$) earns barley a positive return on the extra volume. Barley is therefore assigned a bid of $0.114$ CAD per bushel, which reduces its loss from $2.16$ to $1.90$ million CAD.

Table~\ref{tab:comparison} compares the two cases. Without regulation, the carrier ships the three commodities with positive margins: crude oil ($115.4$ million bushels), corn ($396.8$), and barley ($137.8$), which together fill the network. Oat is not shipped at all, and social welfare is $65.0$ million CAD. This equals the welfare of the null policy $P_1$ on the Pareto frontier (Table~\ref{tab:pareto-summary}), since $P_1$ also sets no quota and the robust stage ships barley at its maximum bid.

Under regulation, the oat quota displaces all of the barley and $110.2$ million bushels of corn, delivers $248.0$ million bushels of oat to the market, and raises social welfare to $288.2$ million CAD.

\begin{table}[H]
\centering
\caption{Shipped quantities, welfare accounting, and player payoffs with and without regulation.}
\label{tab:comparison}
\renewcommand{\arraystretch}{1.25}
\begin{tabularx}{\textwidth}{l Y Y}
\toprule
\textbf{Quantity} & \textbf{Unregulated} & \textbf{Regulated (EGT)} \\
\midrule
Crude oil shipped (M bu)   & $115.4$ & $115.4$ \\
Corn shipped (M bu)        & $396.8$ & $286.6$ \\
Barley shipped (M bu)      & $137.8$ & $0.0$   \\
Oat shipped (M bu)         & $0.0$   & $248.0$ \\
Total shipped (M bu)       & $650.0$ & $650.0$ \\
Network utilization        & $100.0\%$ & $100.0\%$ \\
\midrule
Social welfare $f_1 = \sum_i \beta_i u_i$ (M CAD)  & $65.0$   & $288.2$ \\
Enforcement cost $f_2 = \sum_i K_i(L_i)$ (M CAD)   & $0.0$    & $4.8$   \\
Net regulator utility $f_1 - f_2$ (M CAD)          & $65.0$   & $283.5$ \\
\midrule
Carrier profit (M CAD)          & $57.7$      & $28.2$      \\
Total supplier profit (M CAD)   & $1{,}555.4$ & $1{,}580.8$ \\
\bottomrule
\end{tabularx}
\end{table}

\begin{figure}[H]
\centering
\includegraphics[width=\textwidth]{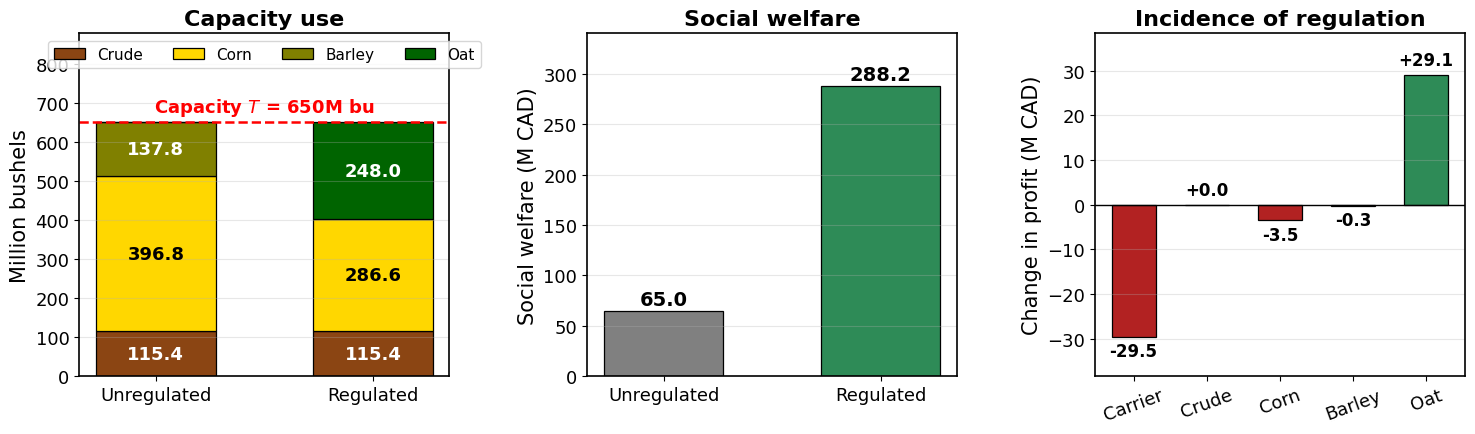}
\caption{Regulated versus unregulated market. \\Left: capacity use by commodity; both markets fill the network. \\Center: social welfare. \\Right: the incidence of regulation, measured as the change in each party's profit relative to the unregulated benchmark.}
\label{fig:comparison}
\end{figure}

Quota regulation raises social welfare from $65.0$ million CAD to $288.2$ million CAD, a gain of $223.2$ million CAD, at an enforcement cost of $4.8$ million CAD, about 2\% of the gain. Both markets use the full capacity of the network, so the gain comes entirely from what is shipped. Without regulation, the capacity left after crude oil and corn goes to barley, and oat is excluded.

In terms of economic incidence of the regulation, there are interesting conclusions to be drawn. The carrier's profit decreases from $57.7$ million CAD to $28.2$ million CAD since the carrier is forced to carry oat, on which its margin is negative, in place of corn and barley, which means a decrease of $29.5$ million CAD that is the private cost the carrier incurs in case of regulation. The total supplier profit increases by $25.4$ million CAD, as the oat supplier, completely excluded in the unregulated market, gains the surplus that emerges due to the guaranteed transport. It means that in this case the regulation results in redistribution of the surplus from the carrier to the previously excluded supplier and generates the net social gain not taken into account by private stakeholders. This is the case described by the political economy of the 2013-2014 Canadian rail crisis, where carriers objected to forced grain shipments due to their negative impact on more profitable freight. In addition, it is the case when intervention becomes necessary since some commodity of a high social value is going to be excluded from the market due to purely profit considerations.

Barley loses its access to the network once the oat quota is in place. Corn and barley then compete for the remaining $286.6$ million bushels, and corn can always outbid barley because its maximum margin ($0.0342$ CAD per bushel) is about twelve times that of barley ($0.0029$). A barley quota would not raise welfare in this setting. Each bushel of barley shipped under a mandate displaces a bushel of corn with the same social weight ($\beta = 0.1$), so welfare stays the same while enforcement cost rises, and the robust stage therefore assigns barley no quota. The exclusion is a loss to the barley supplier, which pays $2.2$ million CAD in inventory penalties, though not a loss in social welfare as the model measures it. Barley receives a quota only when transport costs rise enough to make its maximum margin negative. Raising its social weight above that of corn does not lead to protection, as Section~\ref{subsec:sensitivity} shows.

\subsection{Sensitivity Analysis}
\label{subsec:sensitivity}

The robustness of the equilibrium is examined by varying six parameters around their baseline values: the social weights of oat and barley, the enforcement-cost scale, the worst-case demand factor, the rail capacity, and the carrier's transport cost. In each case the full framework, the ARO frontier followed by the evolutionary dynamics, is re-solved. Unless otherwise stated, all parameters other than the one under study remain at their baseline values.

\begin{figure}[H]
\centering
\includegraphics[width=\textwidth]{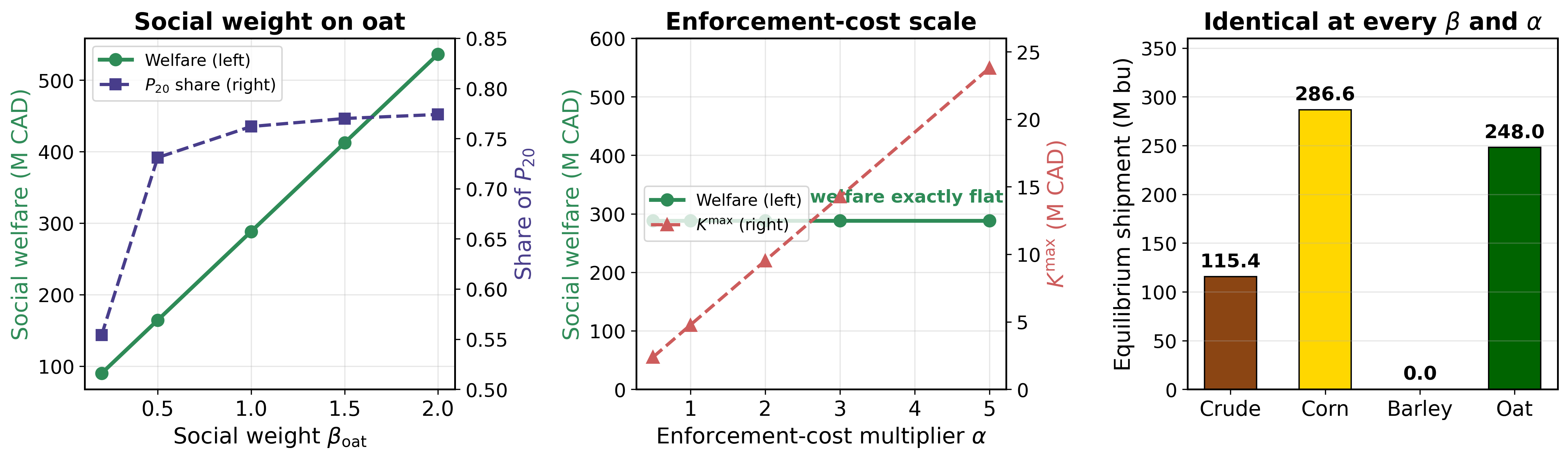}
\caption{Sensitivity to the social weight on oat and to the enforcement-cost multiplier. \\Left: welfare rises linearly in $\beta_{\mathrm{oat}}$ while the share held by $P_{20}$ increases. \\Center: welfare is invariant in $\alpha$, with only the maximum feasible budget $K^{\max}$ scaling. \\Right: the equilibrium allocation, identical at every value of both parameters.}
\label{fig:sa-weights}
\end{figure}

\paragraph{Social weight on oat.}
Varying $\beta_{\mathrm{oat}}$ from $0.2$ to $2.0$ leaves the equilibrium allocation entirely unchanged at $115.4$, $286.6$, $0.0$, and $248.0$ million bushels of crude oil, corn, barley, and oat. Only welfare changes were realized, linearly scaling from $89.8$ to $536.2$ million CAD. The reason is that oat is supplied exclusively through the mandated quota and the marginal welfare gain from protecting it exceeds the marginal enforcement cost at every level considered. Even at the lowest weight considered ($\beta_{\mathrm{oat}} = 0.2$, twice that of crude oil or corn), a bushel of oat quota produces more social value than it costs to enforce. What changes is how strongly the regulator population concentrates on $P_{20}$, whose share rises from $0.554$ at $\beta_{\mathrm{oat}} = 0.2$ to $0.774$ at $\beta_{\mathrm{oat}} = 2.0$. A higher social weight leaves the selected policy unchanged and makes its selection stronger (Figure~\ref{fig:sa-weights}).

\paragraph{Enforcement-cost scale.}
The multiplier $\alpha$ adjusts the weight of the linear and convex portions of $K_i(L_i)$. Increasing $\alpha$ from $0.5$ to $5.0$ leaves the equilibrium unchanged because the regulator keeps $P_{20}$, the shipments profile stays the same, and the welfare equals $288.2$ million CAD for all $\alpha$ values considered. Only the budget axis is scaled, resulting in a new value of $K^{\max}$ between $2.38$ and $23.79$ million CAD, respectively. Invariance holds structurally and not numerically because $K^{\max}$ represents the amount of spending that leads to full protection of oat. Proportionally rescaling the cost function scales the expenditure and budget in such a way that full protection of oat remains feasible at each value of the multiplier.

\begin{figure}[H]
\centering
\includegraphics[width=\textwidth]{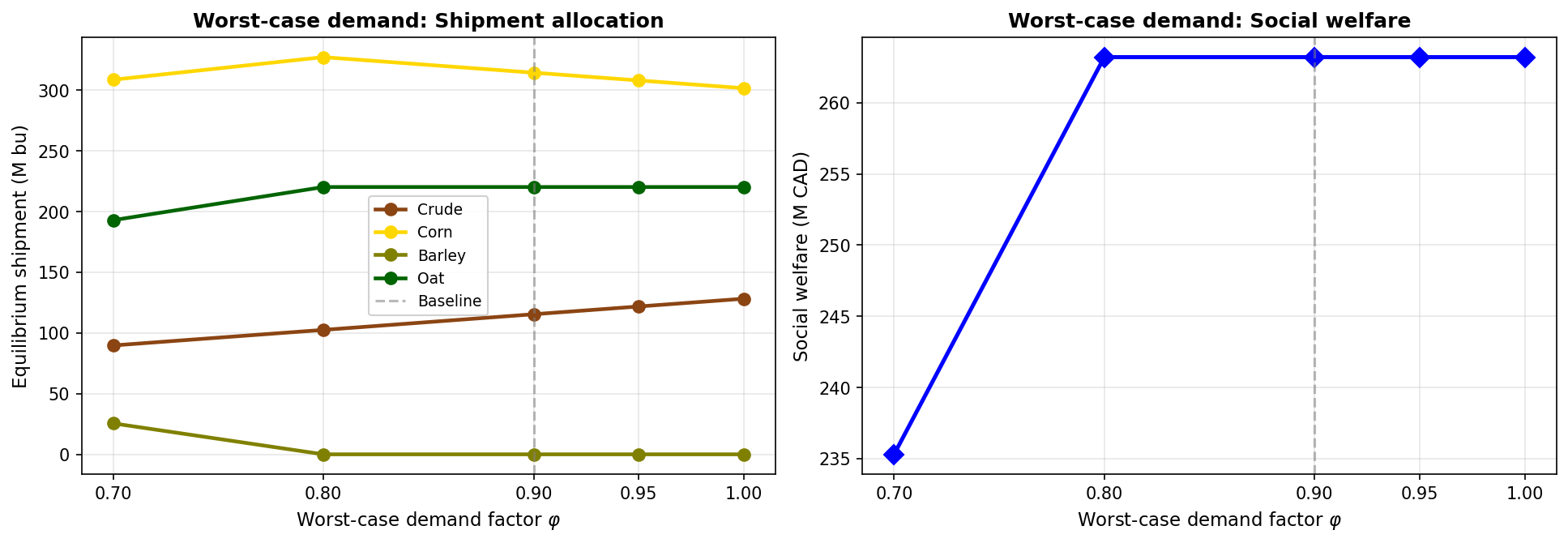}
\caption{Sensitivity to the worst-case demand factor $\varphi$. \\Left: equilibrium shipment of each commodity. \\Right: social welfare.}
\label{fig:sa-demand}
\end{figure}

\paragraph{Worst-case demand factor.}
The factor $\varphi$ determines how conservatively the regulator provisions against uncertainty in commodity availability. Varying it from $0.70$ to $1.00$ shifts the welfare from $232.7$ to $313.0$ million CAD (Table~\ref{tab:sa-phi} and Figure~\ref{fig:sa-demand}). The structure of the equilibrium is stable throughout. The regulator retains $P_{20}$ at a share between $0.759$ and $0.792$, barley never enters the allocation, and oat is protected at exactly its worst-case availability $\varphi D_{\mathrm{oat}}$, which increases from $192.9$ to $275.6$ million bushels.

\begin{table}[H]
\centering
\footnotesize
\caption{Sensitivity to the worst-case demand factor $\varphi$. Oat is protected at its full worst-case availability throughout; under severe uncertainty capacity ceases to bind and part of the network remains idle.}
\label{tab:sa-phi}
\renewcommand{\arraystretch}{1.25}
\setlength{\tabcolsep}{5pt}
\begin{tabularx}{\textwidth}{l Y Y Y Y Y Y}
\toprule
\textbf{$\varphi$} & \textbf{Crude (M bu)} & \textbf{Corn (M bu)} 
& \textbf{Barley (M bu)} & \textbf{Oat (M bu)} & \textbf{Idle (M bu)} 
& \textbf{Welfare (M CAD)} \\
\midrule
$0.70$ & $89.7$  & $308.6$ & $0.0$ & $192.9$ & $58.7$ & $232.7$ \\
$0.80$ & $102.5$ & $327.0$ & $0.0$ & $220.5$ & $0.0$  & $263.4$ \\
$0.90$ & $115.4$ & $286.6$ & $0.0$ & $248.0$ & $0.0$  & $288.2$ \\
$0.95$ & $121.8$ & $266.4$ & $0.0$ & $261.8$ & $0.0$  & $300.6$ \\
$1.00$ & $128.2$ & $246.2$ & $0.0$ & $275.6$ & $0.0$  & $313.0$ \\
\bottomrule
\end{tabularx}
\end{table}

Severe uncertainty produces a distinct effect. At $\varphi = 0.70$ the availability of all commodities contracts by 30\% and the oat quota contracts with it. Because barley still cannot be shipped voluntarily and oat cannot absorb more than its reduced ceiling, $58.7$ million bushels of capacity stand idle and the capacity constraint ceases to bind. Robust provisioning against a severe downside therefore does not lead the regulator to broaden protection across commodities. Provisioning leads to a smaller protected volume and an underused network. The binding constraint changes from capacity to availability, and the welfare loss of $55.5$ million CAD relative to the baseline is the cost of that shift.

\begin{figure}[H]
\centering
\includegraphics[width=\textwidth]{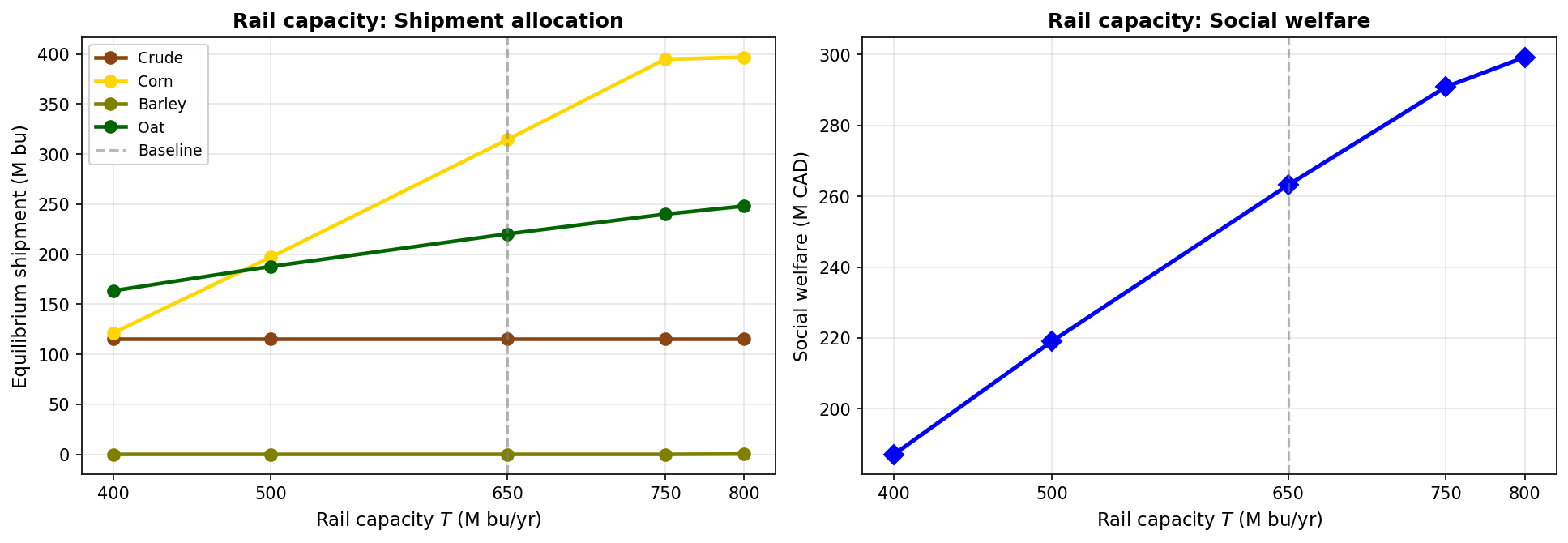}
\caption{Sensitivity to rail capacity $T$. \\Left: crude oil and oat are invariant while corn absorbs the entire adjustment. \\Right: social welfare rises linearly until corn reaches its availability ceiling.}
\label{fig:sa-capacity}
\end{figure}

\paragraph{Rail capacity.}
$T$ was varied between $400$ and $800$ million bushels per year, taking into account the possibilities of disruptions and expansions (Table~\ref{tab:sa-capacity} and Figure~\ref{fig:sa-capacity}). In this range, the regulator always chooses P20, resulting in an evident preference order.

\begin{table}[H]
\centering
\footnotesize
\caption{Sensitivity to rail capacity $T$. Crude oil and oat are unaffected at every level, corn absorbs the entire adjustment.}
\label{tab:sa-capacity}
\renewcommand{\arraystretch}{1.25}
\setlength{\tabcolsep}{5pt}
\begin{tabularx}{\textwidth}{l Y Y Y Y Y Y}
\toprule
\textbf{$T$ (M bu)} & \textbf{Crude (M bu)} & \textbf{Corn (M bu)} 
& \textbf{Barley (M bu)} & \textbf{Oat (M bu)} & \textbf{Idle (M bu)} 
& \textbf{Welfare (M CAD)} \\
\midrule
$400$ & $115.4$ & $36.6$  & $0.0$ & $248.0$ & $0.0$  & $263.2$ \\
$500$ & $115.4$ & $136.6$ & $0.0$ & $248.0$ & $0.0$  & $273.2$ \\
$650$ & $115.4$ & $286.6$ & $0.0$ & $248.0$ & $0.0$  & $288.2$ \\
$750$ & $115.4$ & $386.6$ & $0.0$ & $248.0$ & $0.0$  & $298.2$ \\
$800$ & $115.4$ & $396.8$ & $0.0$ & $248.0$ & $39.8$ & $299.2$ \\
\bottomrule
\end{tabularx}
\end{table}

Crude oil ships in full at every capacity level, including the most severe disruption, because its dominant transport margin guaranties absolute priority. Oat is also invariant at $248.0$ million bushels, since the regulator protects its full worst-case availability regardless of how much capacity is available. Corn absorbs the entire adjustment, falling from $396.8$ million bushels under expansion to $36.6$ million at $T = 400$, a loss of 91\%. The burden of capacity contraction therefore falls neither on the highest-margin commodity nor on the protected one, but entirely on the unprotected commodity that the carrier ships voluntarily.

The marginal social value of capacity is $0.100$ CAD per bushel across the range $T = 400$ to $T = 750$. This is exactly $\beta_{\mathrm{corn}}$. Since corn is the commodity that absorbs the marginal bushel, the capacity is precisely worth its social weight. Above $T = 750$ corn reaches its availability limit of $396.8$ million bushels, the capacity is no longer binding, $39.8$ million bushels are idle and the marginal value drops to $0.020$ CAD per bushel. The shadow price of capacity is therefore the social weight of whichever commodity absorbs the margin, and it collapses once that commodity saturates.

\begin{figure}[H]
\centering
\includegraphics[width=\textwidth]{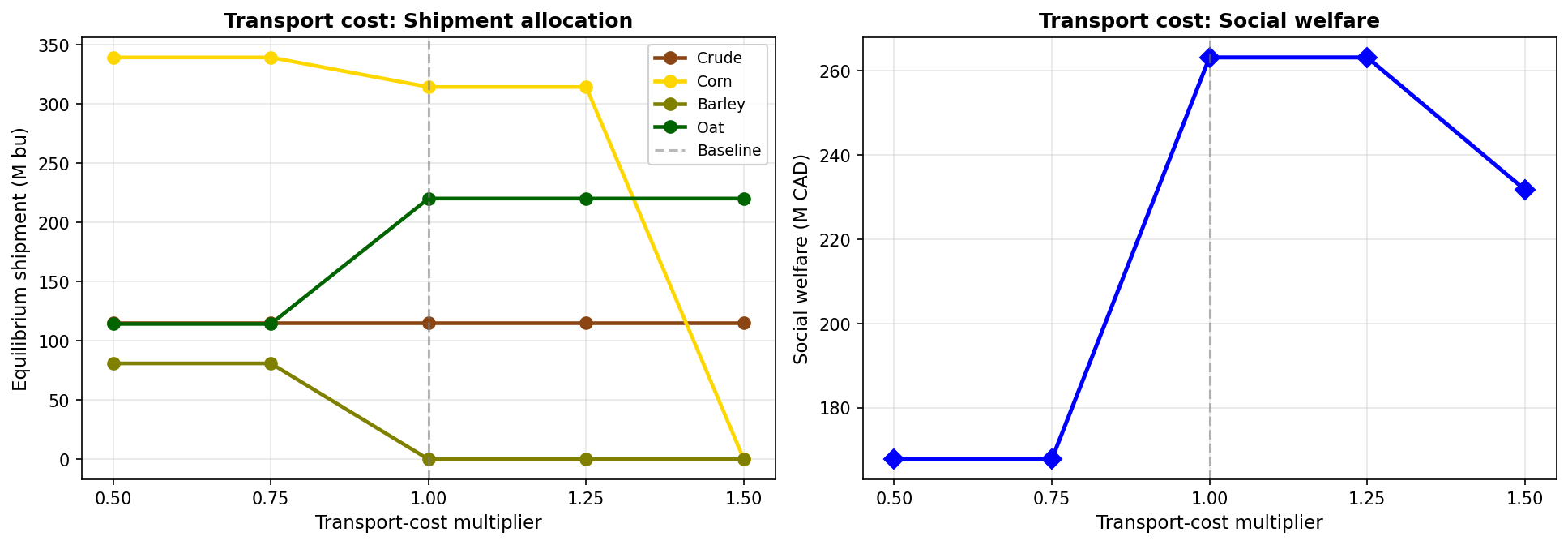}
\caption{Sensitivity to the transport-cost multiplier. \\Left: equilibrium shipment of each commodity. Oat disappears entirely at low cost. \\Right: social welfare is reduced when transport costs fall.}
\label{fig:sa-cost}
\end{figure}

\paragraph{Carrier transport cost.}
The final parameter is the carrier's per-unit operating cost, scaled by a common multiplier to represent a system-wide shock such as a change in fuel or labor costs. Unlike the other parameters, the transport cost changes the signs of the carrier's margins, and hence the commodities that require protection at all. This produces three regimes and a strongly non-monotone welfare response (Table~\ref{tab:sa-cost} and Figure~\ref{fig:sa-cost}).

\begin{table}[H]
\centering
\footnotesize
\caption{Sensitivity to the carrier transport-cost multiplier. At low cost no commodity requires protection and regulation ceases. At high cost both barley and oat require quotas. Under the high-cost multipliers the split of residual capacity between the equally weighted commodities is not uniquely determined.}
\label{tab:sa-cost}
\renewcommand{\arraystretch}{1.25}
\setlength{\tabcolsep}{5pt}
\begin{tabularx}{\textwidth}{l Y Y Y Y Y Y}
\toprule
\textbf{Cost mult.} & \textbf{Oat quota (M bu)} & \textbf{$K^{\max}$ (M CAD)} 
& \textbf{Policy} & \textbf{Share} & \textbf{Carrier profit (M CAD)} 
& \textbf{Welfare (M CAD)} \\
\midrule
$0.50$ & $0.0$   & $0.00$  & $P_{1}$  & $0.050$ & $153.8$  & $65.0$  \\
$0.75$ & $0.0$   & $0.00$  & $P_{1}$  & $0.050$ & $103.7$  & $65.0$  \\
$1.00$ & $248.0$ & $4.76$  & $P_{20}$ & $0.762$ & $28.2$   & $288.2$ \\
$1.25$ & $248.0$ & $12.37$ & $P_{20}$ & $0.196$ & $-52.8$  & $288.2$ \\
$1.50$ & $248.0$ & $12.37$ & $P_{20}$ & $0.195$ & $101.7$  & $288.2$ \\
\bottomrule
\end{tabularx}
\end{table}

At the baseline cost and higher values, oat is fully protected and welfare is maintained at $288.2$ million CAD. When the transport cost increases to $1.25$ times the baseline, the maximum-bid margin of barley also becomes negative, so both barley and oat require quotas and $K^{\max}$ rises to $12.37$ million CAD. Welfare remains unchanged because the additional protected volume carries the same social weight as the corn it displaces. However, two other changes are observed. The carrier's profit decreases to $-52.8$ million CAD, as it is compelled to carry two loss-making commodities at once, and the regulator's dominant policy share falls to $0.196$, indicating a diffuse equilibrium in which many quota policies deliver nearly identical utility. Given the multipliers involved, the division of the residual capacity among the equally valued commodities cannot be uniquely determined. Only the quota for oats and welfare can be determined.

A more important result arises when transportation costs decrease. For multipliers smaller than or equal to $0.75$, the reduction in carrier cost is enough to make the maximum-bid margin of oat positive, which is equal to $m_{\mathrm{oat}} = S_{\mathrm{oat}} + g_{\mathrm{oat}} - c_{\mathrm{oat}}$. Therefore, the ARO rule results in the conclusion that none of the goods require protection, which implies $K^{\max} = 0$ and makes the frontier collapse to the null policy without an imposed quota. However, oat is still excluded. In evolutionary equilibrium, oat bids zero, its margin is negative ($-0.0885$), and a profit-maximizing carrier does not ship it. The welfare decreases from $288.2$ to $65.0$ million CAD and the carrier's profit increases to $153.8$ million CAD. The welfare loss is 77\%.

This is a mechanism of regulatory failure rather than a simple efficiency trade-off. Regulation is switched off at the moment the market appears healthiest because the criterion by which the regulator decides whether protection is needed, the margin at the maximum rational bid, diverges from the margin that actually prevails once suppliers bid strategically. A regulator observing falling transport costs and concluding that intervention is no longer warranted would be reasoning correctly from ex ante margins, reaching a flawed outcome.

\paragraph{Social weight on barley.}
The same failure appears when the social weight of barley is raised above that of corn, with all other parameters at baseline (Table~\ref{tab:sa-barley}). Barley never receives a quota. Its margin at the maximum rational bid is positive ($0.0029$ CAD per bushel), so the robust stage treats it as freight the carrier will ship voluntarily and plans to fill the residual capacity with it, displacing corn and nearly all crude oil. At $\beta_{\mathrm{barley}} = 0.5$ this plan promises $447.0$ million CAD of welfare. At the evolutionary equilibrium barley bids below the carrier's cost, corn takes the residual capacity, and barley ships nothing. Realized welfare stays at $288.2$ million CAD for every weight considered, and the shortfall against the planned welfare grows from $39.7$ million CAD at $\beta_{\mathrm{barley}} = 0.2$ to $357.1$ million CAD at $\beta_{\mathrm{barley}} = 1.0$. A higher social value for barley changes what the regulator expects to achieve but not what the market delivers, since the screening criterion classifies barley as commercially viable.

\begin{table}[H]
\centering
\footnotesize
\caption{Sensitivity to the social weight on barley. Planned values are those of the dominant frontier policy; realized values are those of the evolutionary equilibrium.}
\label{tab:sa-barley}
\renewcommand{\arraystretch}{1.25}
\setlength{\tabcolsep}{5pt}
\begin{tabularx}{\textwidth}{l Y Y Y Y Y Y}
\toprule
\textbf{$\beta_{\mathrm{barley}}$} & \textbf{Barley quota (M bu)} & \textbf{Barley planned (M bu)} & \textbf{Welfare planned (M CAD)} & \textbf{Barley shipped (M bu)} & \textbf{Welfare realized (M CAD)} & \textbf{Shortfall (M CAD)} \\
\midrule
$0.1$ & $0.0$ & --- & $288.2$ & $0.0$ & $288.2$ & $0.0$   \\
$0.2$ & $0.0$ & $396.8$ & $327.9$ & $0.0$ & $288.2$ & $39.7$  \\
$0.5$ & $0.0$ & $396.8$ & $447.0$ & $0.0$ & $288.2$ & $158.7$ \\
$1.0$ & $0.0$ & $396.8$ & $645.4$ & $0.0$ & $288.2$ & $357.1$ \\
\bottomrule
\end{tabularx}
\end{table}

\subsection{Managerial Insights}\label{subsec:managerial}

The results yield five insights of direct relevance to regulators, rail carriers, and commodity producers.

\begin{enumerate}[leftmargin=*, itemsep=2pt]
    \item \textbf{Target intervention narrowly, but fully.} Quota protection delivers large welfare gains only for commodities that are simultaneously socially critical and commercially unviable to transport. Moreover, when such a commodity exists, the analysis favors protecting its entire available volume. In our case study, the marginal enforcement dollar returns 28 times its cost even under full protection.

    \item \textbf{Judge the need for protection on realized bids, not ex ante margins.} A commodity may appear commercially viable at the highest bid the supplier could rationally post and still be excluded at the bid that actually emerges. Screening on ex ante margins in our case study caused the regulation to withdraw entirely when transport costs fell, reducing social welfare by 77\%. Raising the social weight of barley above that of corn exposed the same gap: barley was never protected, and the welfare planned for it did not materialize.

    \item \textbf{Do not be deterred by administrative cost.} In the case studied, the equilibrium was invariant across a tenfold range of enforcement costs. The cost of enforcement determines what protection is offered and not how much is provided.

    \item \textbf{Expect capacity contraction to fall on unprotected voluntary freight.} When disruption occurred in the studied case, high-margin and quota-protected commodities were both unaffected and the entire adjustment fell on corn, which lost 91\% of its volume at the lowest capacity. The  marginal social value of the capacity is equal to the social weight of the commodity that absorbs it.

    \item \textbf{Anticipate that severe supply uncertainty shrinks protection rather than broadening it.} In our case study, conservative provisioning reduced the protected volume in step with availability and left capacity idle. The binding constraint shifted from capacity to availability, and the policy response was a smaller mandate rather than a more diversified one.
\end{enumerate}

\section{Conclusions}\label{sec:conclusion}

This paper developed an integrated framework for the design of minimum shipment quota regulation in rail freight markets under uncertainty. The framework combines adaptive robust optimization, multi-objective optimization, and evolutionary game theory, where the robust optimization stage generates a Pareto frontier of quota policies that trade social welfare against enforcement cost while remaining robust to uncertain commodity availability, and the evolutionary stage embeds this frontier as the regulator's strategy space in a three-population game in which the regulator, the suppliers, and the carrier co-evolve under replicator dynamics. To the best of our knowledge, this is the first model to bring these three methodologies together in a transportation regulation setting, and the first to select a quota policy from a robust Pareto frontier through strategic adaptation of the regulated market. The framework was calibrated to a Canadian rail case study covering crude oil, corn, barley, and oat, and its behavior was examined through a systematic sensitivity analysis.

The results establish that quota regulation can produce a substantial net social benefit. In the Canadian case, regulation raised social welfare from $65.0$ to $288.2$ million CAD, a gain of $223.2$ million CAD, at an enforcement cost of $4.8$ million CAD, about 2\% of the gain. The gain arises because oat combines high social value with a negative transport margin, so a profit-driven market excludes it entirely and it reaches its destination only through a mandated quota. Regulation shifts surplus from the carrier, whose profit falls from $57.7$ to $28.2$ million CAD, to the previously excluded oat supplier, and produces a net social gain that neither private party takes into account.

From the sensitivity analysis, it follows that the optimal solution shows no sensitivity to either the social weight applied to the protected commodity or to the extent of the enforcement cost. However, it exhibits high sensitivity to physical capacity, to the depth of demand uncertainty, and to the carrier’s cost structure. There are three findings that deserve special attention. First, capacity serves as the primary instrument; the shortage does not unduly penalize the highest-value commodity or the protected one, but rather the unprotected commodity, which the carrier transports at its own discretion. That is, the volume of corn is reduced by 91\% due to an extreme capacity shortage, just as was the case with the motivation for the Fair Rail for Grain Farmers Act. The marginal social value of the capacity is equal to the social weight of the commodity that bears the marginal bushel. Second, and most importantly, a decrease in the carrier’s transportation cost can lead to deregulation even where regulation would have been needed. Provided that the cost decrease is such that the margin of the protected commodity at the maximal rational bid becomes positive, then the robust stage finds that there is no need for protection and does not impose any quota. However, the commodity continues to make a zero bid, the margin is negative, and the commodity is excluded from transportation. Social welfare decreases by 77\%, while the carrier’s profit increases. The same mechanism leaves barley unprotected when its social weight exceeds that of corn; the robust stage expects barley to be shipped voluntarily, and none is shipped at equilibrium. Third, severe demand uncertainty does not lead the regulator to spread protection across commodities. It reduces the protected volume in step with contracted availability, the binding constraint changes from capacity to availability, and part of the network remains idle.

The second finding has implications beyond the case at hand. It identifies a mechanism of regulatory failure that arises from the gap between the margin a regulator can observe before bidding takes place and the margin that prevails once suppliers adapt. A regulator screening candidates for protection on ex ante margins would reason correctly from the information available and still reach the wrong policy. Only a framework in which bidding behavior is an outcome rather than an assumption exposes this gap, which is the practical case for coupling robust design with adaptive dynamics rather than treating either in isolation.

A look at the research limitations highlights directions for future research. The model represents the rail network by a single aggregate capacity constraint, abstracting from route-level and seasonal congestion. Extending the framework to a spatially disaggregated network will allow the regulator to introduce geographically targeted quota and incorporate equity considerations across regions. Additionally, the uncertainty set was restricted to commodity availability using an interval representation. Richer structures, including correlated shocks to prices, costs, and capacity, and alternative representations such as distributionally robust or data-driven sets \citep{cheng2026distributionally}, will broaden the applicability of the approach. Furthermore, the robust stage screens commodities for protection using the margin at the maximum rational bid, which is the conservative criterion available before bidding occurs. An iterative scheme in which the frontier is recomputed against the bids the dynamics generates will close the gap, at the cost of a considerably harder computational problem. On another note, the case study considered a single carrier and four representative commodities. Generalizing to multiple competing carriers and a larger commodity set will capture horizontal competition among transporters. In addition, quotas in this paper are set for each commodity, whereas the Fair Rail for Grain Farmers Act imposed a single quota on total grain shipments. Under an aggregate quota the carrier chooses which grain to ship, and extending the model in this direction would show how grains with thin transport margins, such as barley, fare under the form of regulation used in practice. The framework was also limited to a single planning period. A multi-period extension incorporating inventory carryover, the time value of money, and perishability will reflect the dynamic nature of harvest and shipment cycles. Market prices are held fixed as well, so a commodity excluded from transport does not become more valuable in its destination market. If prices responded to shipped quantities through the inverse demand of Eq.~\eqref{eq:inverse-demand}, a scarce commodity such as barley could bid more for capacity over successive periods, and part of its exclusion might correct itself. Moreover, when several commodities have the same relative social importance, the allocation of surplus capacity between them is not unique, and another criterion, such as minimization of variance in share distributions, will be required to eliminate the indeterminacy.

\bibliographystyle{elsarticle-harv}
\bibliography{02_References} 

\end{document}